\documentclass[preprint]{elsarticle}

\usepackage{amsmath,amsfonts}
\usepackage{amssymb}

\newcommand{\R}{\mathbb R}

\newtheorem{theorem}{Theorem}
\newtheorem{lemma}{Lemma}
\newtheorem{example}{Example}
\newproof{proof}{Proof}

\usepackage{epsfig}

\newcommand{\cH}{\mathcal{H}}
\newcommand{\cD}{\mathcal{D}}

\newcommand{\cS}{\mathcal{S}}
\newcommand{\cN}{\mathcal{N}}
\newcommand{\bcH}{\overline{\cH}}

\newcommand{\bu}{\overline{u}}

\newcommand{\bcS}{\overline{\cS}}
\newcommand{\bcN}{\overline{\cN}}

\usepackage[dvipsnames]{xcolor}

\begin{document}

\title{Preserving energy resp. dissipation in numerical PDEs using the ``Average Vector Field'' method}

\author[ntnu]{E. Celledoni}\ead{elenac@math.ntnu.no}
\author[karlsruhe]{V. Grimm}\ead{volker.grimm@kit.edu}
\author[ifs]{R. I. McLachlan}\ead{r.mclachlan@massey.ac.nz}
\author[lat]{D. I. McLaren}\ead{d.mclaren@latrobe.edu.au}
\author[lat,irl]{D. O'Neale}\ead{d.oneale@irl.cri.nz}
\author[ntnu]{B. Owren}\ead{brynjulf.owren@math.ntnu.no}
\author[lat]{G. R. W. Quispel}\ead{r.quispel@latrobe.edu.au}

\address[ntnu]{Department of Mathematical Sciences, NTNU
7491 Trondheim, Norway.}
\address[karlsruhe]{Karlsruhe Institute of Technology (KIT), Institut f\"ur Angewandte und Numerische Mathematik, D-76128 Karlsruhe, Germany. }
\address[ifs]{Institute of Fundamental Sciences, Massey University,
Private Bag 11-222\\ Palmerston North, New Zealand.}
\address[lat]{Department of Mathematics and Statistics, La Trobe University,
Victoria 3086, Australia.}
\address[irl]{Industrial Research Ltd, PO Box 31310, Lower Hutt 5040, New Zealand}

\begin{abstract}
We give a systematic method for discretizing Hamiltonian partial differential equations
(PDEs) with constant symplectic structure, while preserving their energy
exactly. The same method, applied to PDEs with constant dissipative structure,
also preserves the correct monotonic decrease of energy. The method is
illustrated by many examples. In the Hamiltonian case these include:  the
sine--Gordon, Korteweg--de Vries, nonlinear Schr\"{o}dinger, (linear)
time-dependent Schr\"{o}dinger, and Maxwell equations. In the dissipative case
the examples are:  the Allen--Cahn, Cahn--Hilliard, Ginzburg--Landau, and heat equations.
\end{abstract}


\maketitle



\section{Introduction} ``The opening line of {\em Anna Karenina}, `All happy families resemble one
another, but each unhappy family is unhappy in its own way', is a useful
metaphor for the relationship between
computational ordinary differential equations (ODEs) and computational partial
differential equations (PDEs). ODEs are a happy family -- perhaps they do not
resemble each other, but, at the very least, we can treat them by a 
relatively small compendium of computational techniques\dots
PDEs are a huge and motley collection of
problems, each unhappy in its own way'' (Quote from A.~Iserles' book
\cite{Iserlesbook}).

Whereas there is much truth in the above quote, in this paper we set out to
convince the reader that, as far as conservation or dissipation of energy is
concerned, many PDEs form part of one big happy family (cf. also \cite{Klainerman08PCP}) that, after a very straightforward and uniform 
semi-discretization, may actually be solved by a single unique geometric 
integration method -- the so-called average vector field method -- while 
 preserving the correct conservation, respectively, dissipation of energy. The concept of `energy' has far-reaching importance throughout the physical sciences \cite{Feyn1}. Therefore a single procedure, as presented here, that correctly conserves, resp. dissipates, energy for linear as well as nonlinear, low-order as well as high-order, PDEs would seem to be worth while.

Energy-preserving schemes have a long history, going back to Courant, Friedrichs, and Lewy's cunning derivation \cite{Courant28} of a discrete energy conservation law for the 5-point finite difference approximation of the wave equation which they used to prove the scheme's convergence. The conservation law structure of many PDEs is considered fundamental to their derivation, their behaviour, and their discretization.
Li and Vu-Quoc \cite{Li95} give a historical survey of energy-preserving methods for PDEs and their applications, especially to nonlinear stability. What is relevant to us here is that many of these methods (e.g. \cite{Celledoni08,Fei95,Fei91,Furihata99,Matsuo07,Matsuo01,Matsuo09,Ringler10,Yaguchi10}) have an ad-hoc character and are not completely systematic either in their derivation or in their applicability; in contrast, the method discussed here (Eq. \eqref{avf} below) is completely systematic, applies to a huge class of conservation and dissipative PDEs, and depends functionally only on the PDE itself, not its energy. In some cases it reduces to previously studied methods, for example, it reproduces one of Li and Vu-Quoc's schemes \cite{Li95} for the nonlinear wave equation. Even in these cases, however, it sheds considerable light on the actual structure of the scheme and the origin of its conservative properties. See also the discussion comparing different constructions of energy-preserving integrators in \cite{Dahlby11}.

We consider evolutionary PDEs with independent variables 
$(x,t) \in \R^d \times \R$, functions $u$
belonging to a Hilbert space $\mathcal{B}$ with values\footnote{Although it is generally real-valued, the function $u$ may also be complex-valued, for example, the nonlinear Schr\"{o}dinger equation.} 
$u(x,t) \in \R^m$, and PDEs of the form
\begin{equation} \label{basicpde}
\dot{u}=\cD  \frac{\delta \cH}{\delta u},
\end{equation}
where $\cD$ is a constant linear differential operator, the dot denotes 
$\frac{\partial}{\partial t}$, and 
\begin{equation}
\cH [u] = \int_{\Omega}^{} H(x;u^{(n)}) \; dx
\end{equation}
where $\Omega$ is a subset of $\R^d \times \R$, and $dx = dx_1dx_2 \dots dx_d$.
 $\frac{\delta \cH }{\delta u}$ is the 
variational derivative of $\cH$ in the
sense that
\begin{equation}
   \frac{d}{d \epsilon} \cH [u+\epsilon v] 
   \big|_{\epsilon=0}
   =
   \left<\frac{\delta \cH}{\delta u}, v\right>
\end{equation}
for all $u,v \in \mathcal{B}$ (cf. \cite{Olver93}), where $\left<,\right>$ is the inner product in $\mathcal{B}$. For example, if 
$d=m=1$, $\mathcal{B}=L^2(\Omega)$, and
\begin{equation}
  \cH[u] =\int_{\Omega} H(x;u,u_x,u_{xx},\dots) \, dx,
\end{equation}
then
\begin{equation}
\frac{\delta \cH}{\delta u} = \frac{\partial H}{\partial u} - 
\partial_x \left( \frac{\partial H}{\partial u_x} \right) + \partial_x^2 
\left( \frac{\partial H}{\partial u_{xx}} \right) - \cdots, 
\end{equation}
when the boundary terms are zero.

Similarly, for general $d$ and $m$, we obtain
\begin{equation} \label{vardivmain}
\frac{\delta \cH}{\delta u_l} = \frac{\partial H}{\partial u_l} - 
\sum_{k=1}^{d}
\frac{\partial}{\partial {x_k}} \left( \frac{\partial H}{\partial u_{l,k}} \right) + \dots, \quad l=1, \dots,m. 
\end{equation}

We consider Hamiltonian systems of the form 
(\ref{basicpde}), where $\cD$ is a constant skew symmetric operator 
(cf. \cite{Olver93}) and $\cH$ the energy (Hamiltonian). In this case,
we prefer to designate the differential operator in (\ref{basicpde}) with 
$\cS$ instead of $\cD$. The PDE 
preserves the energy because $\cS$ is skew-adjoint with respect to the 
$L^2$ inner product, i.e.
\begin{equation}
\int_{\Omega} u\cS u \, dx = 0, \quad \forall u \in \mathcal{B}.
\end{equation}
The system (\ref{basicpde}) has $\mathcal{I}:\mathcal{B} \rightarrow \R$ as an
integral if $\dot{\mathcal{I}}=\int_{\Omega} \frac{\delta \mathcal{I}}{\delta u}
\cS \frac{\delta \cH} {\delta u} \, dx = 0$. 

\noindent Integrals
$\mathcal{C}$ with $\cD \frac{\delta \mathcal{C}}{\delta u}=0$ are
called Casimirs.

Besides PDEs of type (\ref{basicpde}) where $\cD$ is skew-adjoint, we 
also consider PDEs of type (\ref{basicpde}) where $\cD$ is a constant negative (semi)definite operator  
with respect to the $L^2$ inner product, i.e.
\begin{equation}
\int_{\Omega} u \cD u \, dx \leq 0, \quad  \forall u \in \mathcal{B}.
\end{equation}
In this case, we prefer to designate the differential operator $\cD$ with $\cN$ 
and the function $\cH$ is a Lyapunov function, 
since then the system (\ref{basicpde}), i.e.  
\begin{equation} \label{lyappde}
\dot{u}=\cN \frac{\delta \cH}{\delta u},
\end{equation}
has $\cH$ as a Lyapunov function, i.e. 
$\dot{\cH}=\int_{\Omega} \frac{\delta \cH}{\delta u}
\cN \frac{\delta \cH}{\delta u} \, dx \leq 0.$
We will refer to systems (\ref{basicpde}) with a skew-adjoint $\cS$ and an
energy $\cH$ as conservative and to systems (\ref{basicpde}) with a negative (semi)definite operator 
$\cN$ and a Lyapunov function $\cH$ as dissipative. Note that the operator $\cN$ need {\it not} be self-adjoint. (In Example \ref{ex:gl}, the Ginzburg--Landau equation, $\cN = \partial_x + \epsilon\partial_x^2$ is not self-adjoint.)

Conservative PDEs (\ref{basicpde}) can be semi-discretized in 
``skew-gradient'' form
\begin{equation} \label{basicode}
  \dot{u}=\bcS \nabla \bcH (u), \qquad \bcS^T=-\bcS,\qquad u\in{\mathbb R}^k
\end{equation}
when $\cD = \cS$ is skew-adjoint. Here, and in the 
following, we will always denote the discretizations with bars.
$\bcH\colon{\mathbb R}^k\to{\mathbb R}$ is chosen in such a way that  $\bcH \Delta x$ is an approximation to $\cH$. 
\medskip
\begin{lemma}
Let
\begin{equation} 
 \cH [u] = \int_{\Omega}^{} H(x ; u^{(n)} ) dx,
\end{equation}
and let $\bcH \Delta x$ be any consistent (finite difference) approximation to $\cH$ (where $\Delta x := \Delta x_1\Delta x_2 \dots \Delta x_d$) with $N$ degrees of freedom. Then in the finite-dimensional Hilbert space 
$\R^N$ with the Euclidean inner product,  the   variational derivative $\dfrac{\delta}{\delta u}(\bcH\Delta x)$ is given by $\nabla \bcH$.
\end{lemma}
\medskip \noindent
\begin{proof}
We denote the consistent (finite difference) discretization of
$
\cH [u] = \int_{\Omega}^{} H(x;u^{(n)}) \; dx
$
by
$
\bcH' (\bu)= \bcH(\bu) \Delta x $ where $\bu\in\R^N $ denotes the discrete values of
$u$, in the multidimensional case after choosing an ordering.
The variational derivative $\dfrac{\delta
\cH}{\delta u}$ is then given by
\begin{eqnarray} \label{vardivnew}
\frac{d}{d \varepsilon} \bcH'(\bu + \varepsilon \overline{v})
\big|_{\varepsilon = 0} &=& \sum_{n=1}^{N} \left( \dfrac{\delta
\bcH'}{\delta \bu} \right)_n \cdot \overline{v}_n \; \Delta x \\
& = & \left( \overline{v} \cdot \dfrac{\delta \bcH'}{\delta \bu} \right)
\Delta x .  \notag
\end{eqnarray}
It holds for the directional derivative that
\begin{equation} \label{vdpfnew}
\frac{d}{d \varepsilon} \bcH'(\bu + \varepsilon \overline{v})
\big|_{\varepsilon = 0} = (\overline{v} \cdot \nabla \bcH').
\end{equation}
Since Eqs (\ref{vardivnew}) and (\ref{vdpfnew}) must hold for all vectors
$\overline{v}$, we have
\begin{equation}
\dfrac{\delta \bcH'}{\delta \bu} \Delta x \equiv \nabla \bcH' ,
\qquad \mbox{and hence} \qquad
\dfrac{\delta \bcH'}{\delta \bu}  \equiv \nabla \bcH .
\end{equation}
\end{proof}

It is worth noting that the above lemma also applies directly when the approximation to $\cH$ is obtained by a spectral discretization, since such an approximation can be viewed as a finite difference approximation where the finite difference stencil has the same number of entries as the number of grid points on which it is defined.

The operator $\nabla$ is the standard gradient, which replaces the variational derivative because we are now working in a finite (although large) number of dimensions (cf. e.g. (\ref{vardivmain})). 
When dealing with (semi-)discrete systems we use the notation $u_{j,n}$ where the index $j$ corresponds to increments in space and $n$ to increments in time. That is, the point $u_{j,n}$ is the discrete equivalent of $u(a+j\Delta x,t_0+n\Delta t)$ where $x\in[a,b]$ and where $t_0$ is the initial time. In most of the equations we present, one of the indices is held constant, in which case, for simplicity, we drop it from the notation. For example, we use $u_j$ to refer to the values of $u$ at different points in space and at a fixed time level.

\begin{theorem}
Let $\bcS$ (resp. $\bcN$) be any consistent constant skew (resp. negative-definite) matrix approximation to $\cS$ (resp. $\cN$). Let $\bcH \Delta x$ be any consistent (finite difference) approximation to $\cH$. Finally, let 
\begin{equation} \label{discf}
f(u) := \bcS \nabla \bcH (u) \qquad (\mbox{resp. } f(u) := \bcN \nabla \bcH (u) ),
\end{equation}
and let $u_n$ be the solution of the average vector field (AVF) method
\begin{equation} \label{avf}
\frac{u_{n+1}-u_n}{\Delta t}= \int_0^1  f((1-\xi)u_n+\xi u_{n+1}) \, d\xi,
\end{equation}
applied to equation (\ref{discf}).  Then the semidiscrete energy $\bcH$ is preserved exactly (resp. dissipated monotonically):
\begin{equation} \notag
\bcH (u_{n+1}) = \bcH (u_n) \qquad ( \mbox{resp. }  \bcH (u_{n+1}) \leq \bcH (u_n) ).
\end{equation}
\end{theorem}

$\bcH$  is preserved by the flow of $\dot u = \bcS\nabla \bcH(u)$ since
\begin{equation}
\dot{\bcH}=\left( \nabla \bcH \right)^T \bcS \nabla \bcH = 0.
\end{equation}
Discretisations of this type can be given for pseudospectral, finite-element,
Galerkin and finite-difference methods 
(cf. \cite{McLachlan_Robidoux00,MR1870274}); for
simplicity's sake, we will concentrate on finite-difference methods, though we include two examples of pseudospectral methods for good measure. We consider only uniform grids; see \cite{Kitson03} for finite difference discretizations of this type on nonuniform grids, which require inner products
$\left<\overline{u}_n,\overline{v}_n\right> = \overline{u}_n\cdot(M \overline{v}_n)$ with $M\ne I$.

The AVF method was recently \cite{Quispel_McLaren08} shown to preserve the energy $\bcH$ exactly for any vector field $f$ of the form $f(u) = \bcS \nabla \bcH (u)$, where $\bcH$ is an arbitrary function, and $\bcS$ is any \textbf{constant} skew matrix\footnote{The relationship of (\ref{avf}) to Runge-Kutta methods was explored in \cite{Celledonietal08}.}. The AVF method is related to  discrete gradient
methods (cf. \cite{McLachlanetal99}); amongst discrete gradient methods it is distinguished by its features of linear covariance, automatic preservation of linear symmetries, reversibility with respect to linear reversing symmetries, and often by its simplicity.
It is one member of the family of Galerkin methods introduced
by Betsch and Steinmann \cite{Betsch00} for the case of canonical classical mechanical systems, although the dependence of the method on $f$ alone (and not $\cH$) was not realized there.

We remark that although \eqref{avf} does not specify the method in a completely closed form (because of the integral), the same is true for implicit methods, whose description needs to be supplemented by an iterative solver in any implementation. Despite this, broadening the class of methods from explicit to implicit is essential for many differential equations and allows new properties such as A-stability and symplecticity. Similarly, a further broadening to include integrals of the vector field allows the new properties of energy conservation and dissipation.

In many cases---for all linear, polynomial, or scalar terms in $f$, which includes all the examples in this paper---the integral can be evaluated exactly in closed form. In other cases, it can be approximated by quadrature to any desired degree of accuracy.

If $\cD$  is a constant negative-definite operator, then the dissipative PDE (\ref{basicpde})
can be discretized in the form
\begin{equation} \label{baseodelyapi}
\dot{u}=\bcN \nabla \bcH (u), 
\end{equation} 
where $\bcN$ is a negative (semi)definite matrix and $\bcH$ is a discretization  as above.

That is, $\bcH$ is a Lyapunov-function for the semi-discretized system, since
\begin{equation}
\dot{\bcH} = \left( \nabla \bcH \right)^T \bcN \nabla \bcH \leq 0.
\end{equation}
The AVF method (\ref{avf}) again preserves this structure, i.e. we have
\begin{equation}
\bcH (u_{n+1}) \leq \bcH (u_n),
\end{equation}
and $\bcH$ is a Lyapunov function for the discrete system. Taking the scalar
product of (\ref{avf}) with 
$\int_0^1 \nabla \bcH ((1-\xi)u_n+\xi u_{n+1}) \, d\xi$ 
on both sides of the equation yields
\begin{equation}
\frac{1}{\Delta t} \int_0^1 (u_{n+1}-u_n) \cdot\nabla \bcH ((1-\xi)u_n+\xi u_{n+1}) \, d\xi
\leq 0,
\end{equation}
i.e.
\begin{equation}
\frac{1}{\Delta t} \int_0^1 \frac{d}{d \xi} \bcH ((1-\xi)u_n+\xi u_{n+1}) \, d\xi \leq 0,
\end{equation}
and therefore
\begin{equation}
\frac{1}{\Delta t} (\bcH (u_{n+1})-\bcH (u_n))\leq 0.
\end{equation}

Our purpose is to show that the procedure described above, namely
\begin{enumerate}
\item Discretize the energy functional $\cH$ using any (consistent) approximation $\bcH \Delta x$
\item Discretize $\cD$ by a constant skew-symmetric (resp. negative (semi)definite) matrix
\item Apply the AVF method
\end{enumerate}
can be generally applied and leads, in a systematic way, to energy-preserving methods for 
conservative PDEs and energy-dissipating methods for dissipative PDEs. We shall 
demonstrate the procedure by going through several well-known nonlinear and 
linear PDEs step by step. In particular we give examples of how to discretize nonlinear conservative PDEs (in subsection 2.1), linear conservative PDEs (in subsection 2.2), nonlinear dissipative PDEs (in subsection 3.1), and linear dissipative PDEs (in subsection 3.2).

Energy dissipation has been much less studied than energy conservation. Stuart and Humphries \cite{Stuart96} consider conditions for certain linear methods (such as backward Euler) to be dissipative for gradient systems  $\dot x = -\nabla F(x)$, $F(x)\ge 0$, $\lim_{\|x\|\to 0}F(x)=\infty$. This class is however far more restrictive than the dissipative systems \eqref{lyappde}; the restriction from $\cN$ negative (semi)definite and not necessarily self-adjoint to $\cN=-$id alone is substantial \cite{McLachlanetal99}. In particular, backward Euler is not dissipative for systems of the form \eqref{lyappde} (see Example \ref{ex:gl} below).

The method \eqref{avf} is formally second order in time. The relationship between accuracy in space and time is a complicated one that we cannot explore here. Depending on the PDE and the scientific goals the accuracy in time may need to be much less, the same, or much greater than that in space. The temporal order can be increased if necessary by composition \cite{McLachlan02} or by including derivatives of the right hand side \cite{Quispel_McLaren08}.

Steps 1 and 2 yield semidiscretizations that are Hamiltonian (or Poisson) in the conservative case and dissipative in the dissipative case. Thus they could be followed by a symplectic time integrator, so this procedure unifies symplectic and energy-conserving integration of conservative PDEs and allows for systematic comparisons between the two. Other approaches to energy-conserving integration, e.g. \cite{Matsuo07,Matsuo01,Matsuo09}, conflate time and space discretization and obscure the fundamental role played by the semidiscrete energy and its Euclidean gradient.

The choice of a symplectic versus an energy-preserving time integrator is an interesting and delicate question that depends on the PDE and the scientific goals. Some factors favour energy conservation:
\begin{enumerate}
\item
When the (semidiscretized) energy level sets are compact, conserving energy may (via energy stability) improve performance at large time steps; Simo and Gonzalez \cite{SimoGonzalez93} show that the unresolved high frequencies are controlled by exact energy conservation, whereas they can lead to instability in the symplectic midpoint rule. (See also their comparisons in the context of relative equilibria \cite{Gonzalez96} and in nonlinear elastodynamics \cite{SimoTarnow92}.)
\item
If a nonconservative scheme is used for a system of conservation laws and large (e.g.) mass errors result, this is not taken to be a sign that the spatial mesh size should be reduced, or that mass should be globally rescaled back to its original value; rather, it is taken to be a sign that a conservative scheme should have been used. Unitarity in Schr\"odinger equations and orthogonality in spin chains followed a similar course: their preservation was first discovered to be important; then enforced by projection; eventually, intrinsically-conserving schemes were discovered that are now widely used.
\item
In fully chaotic dynamics, symplecticity may be irrelevant. It is not used in statistical mechanics (only energy and volume preservation are required). If the flow on an energy level set is Anosov, then it is structurally stable, even to non-Hamiltonian perturbations.
\item
It is easier to adapt the time step in energy-preserving than in symplectic integration.
\item
It leads to different constraints on the eigenvalues and bifurcations of orbits than symplecticity; for example, energy preservation leads to periodic orbits occurring in 1-parameter families (parameterized by the energy), whereas symplecticity does not.
\end{enumerate}
Some factors favour symplecticity:
\begin{enumerate}
\item Unlike energy, symplecticity is the defining property of Hamiltonian mechanics. It ensures preservation of many phase space features such as genericity of quasiperiodic orbits.
\item It can be essential to capture qualitatively correct long-time limit sets and long-time statistics.
\item In $2n$ dimensions, symplecticity provides $n(2n-1)$ constraints, energy only one.
\item It gives the user a very convenient check on the simulation, namely, the energy error. However, the security this provides may be illusory as neither conservation of energy nor small energy errors in a nonconservative scheme ensure that the simulation is reliable.
\end{enumerate}

\section{Conservative PDEs}
\subsection{Nonlinear conservative PDEs}

\begin{example}
Sine--Gordon equation:

\underline{Continuous:}

\end{example}
\begin{equation} \label{sg:equation}
\frac{\partial^2 \varphi}{\partial t^2} = 
\frac{\partial^2 \varphi}{\partial x^2} - \alpha \sin \varphi. 
\end{equation}
The Sine--Gordon equation is of type (\ref{basicpde}) with
\begin{equation} \label{sgbcH}
\cH = \int 
\left[
  \frac{1}{2} \pi^2 + 
  \frac{1}{2} \left( \frac{\partial \varphi}{\partial x} \right)^2
  + \alpha \left( 1 - \cos \varphi \right)
\right] \, dx,
\end{equation}
where $u := \left( \begin{array}{c} \varphi \\ \pi \end{array} \right)$ and
\begin{equation}
\cS = 
\left(
  \begin{array}{cc}
     0  &  1  \\
    -1  &  0
  \end{array}
\right).
\end{equation}
(Note that it follows that $\pi = \frac{\partial \varphi}{\partial t}$.)

Boundary conditions: periodic, $u(-20,t) = u(20,t)$.

\underline{Semi-discrete: finite differences}\footnote{Summations of the form $\sum_{j}^{}$ mean $\sum_{j=0}^{N-1}$ unless stated otherwise.}
\begin{equation}
\bcH_{fd} = \sum_j
 \left[ \frac{1}{2} \pi_j^2 + \frac{1}{2(\Delta x)^2}(\varphi_{j+1}-\varphi_{j})^2 
 + \alpha \left( 1- \cos \varphi_j \right) \right]. 
\end{equation}
\begin{equation}
\bcS = 
\left(
   \begin{array}{cc}
     0  & \mbox{id} \\
    -\mbox{id} & 0
   \end{array}
\right).
\end{equation}

The resulting system of ordinary differential equations is
\begin{equation}
\left[
\begin{array}{c}
{\bf\dot{\boldsymbol{\varphi}}}\\
{\bf\dot{\boldsymbol{\pi}}}
\end{array}
\right]
=\bcS\nabla\bcH_{fd}=
\left[
\begin{array}{c}
  \boldsymbol{\pi}\\
  \frac{1}{\Delta x^2}L\boldsymbol{\varphi}-\alpha\sin\boldsymbol{\varphi}
\end{array}
\right],
\end{equation}
 where $L$ is the circulant matrix
$$
L=\left[
\begin{array}{cccc}
-2&1&&1\\
1&\ddots&\ddots&\\
&\ddots&\ddots&1\\
1&&1&-2\\
\end{array}
\right].
$$
We have used the bold variables $\boldsymbol{\varphi}$ and $\boldsymbol{\pi}$ for the finite dimensional vectors $[\varphi_1,\varphi_2,\ldots,\varphi_N]^{\top}$, \emph{et cetera}, which replace the functions $\pi$ and $\varphi$ in the (semi-) discrete case. Where necessary, we will write $\boldsymbol{\varphi}_n$, \emph{et cetera} to denote the vector $\boldsymbol{\varphi}$ at time $t_0+n\Delta t$.

The integral in the AVF method can be calculated exactly to give\footnote{For numerical computations, care must be taken to avoid problems when the difference $\boldsymbol{\varphi}_{n+1}-\boldsymbol{\varphi}_{n}$ in the denominator of (\ref{SineGordonAVF}) becomes small. We used the sum-to-product identity $\cos a-\cos b=-2\sin((a+b)/2)\sin((a-b)/2)$ to give a more numerically amenable expression. In Eq. (\ref{SineGordonAVF}), elementwise division of vectors is used.}
\begin{align}
\label{SineGordonAVF}
\frac{1}{\Delta t}&\left[
\begin{array}{c}
  \boldsymbol{\varphi}_{n+1}-\boldsymbol{\varphi}_n\\
  \boldsymbol{\pi}_{n+1}-\boldsymbol{\pi}_n
\end{array}
\right]
=\\
&\left[
\begin{array}{c}
  (\boldsymbol{\pi}_{n+1}+\boldsymbol{\pi}_n)/2\\
  L(\boldsymbol{\varphi}_{n+1}+\boldsymbol{\varphi}_n)/2-\alpha(\cos\boldsymbol{\varphi}_{n+1}-\cos\boldsymbol{\varphi}_n)/(\boldsymbol{\varphi}_{n+1}-\boldsymbol{\varphi}_n)
\end{array}
\right].\notag
\end{align}

\underline{Semi-discrete: spectral discretization}

Instead of using finite differences for the discretization of the spatial derivative in (\ref{sgbcH}), one may use a spectral discretization. This can be thought of as replacing $\varphi$ with its Fourier series, truncated after $N$ terms, where $N$ is the number of spatial intervals, and differentiating the Fourier series. This can be calculated, using the discrete Fourier transform\footnote{In practice, one uses the fast Fourier transform algorithm to calculate the DFTs in $\mathcal{O}(N\log N)$ operations.} (DFT), as $\mathcal{F}_N^{-1}D_N\mathcal{F}_N\boldsymbol{\varphi}$ where $\mathcal{F}_N$ is the matrix of DFT coefficients with entries given by $[\mathcal{F}_N]_{n,k}=\omega_N^{nk},$ $\omega_N=e^{-i2\pi/N}$  Additionally, $[\mathcal{F}_N^{-1}]_{n,k}=\frac{1}{N}\omega_N^{-nk}$ and $D_N$ is a diagonal matrix whose (non-zero) entries are the scaled wave-numbers\footnote{Care must be taken with the ordering of the wave numbers since different computer packages use different effective orderings of the DFT/IDFT matrices in their algorithms. Additionally, one must ensure that all modes of the Fourier spectrum are treated symmetrically --- for $N$ even, this requires replacing the $k=\frac{N}{2}$ entry with zero to give $[0,\ldots,\frac{N}{2}-1,0,\frac{-N}{2}+1,\ldots,-1]$.} $$\text{diag}(D_N)=\frac{2\pi i}{l}\left[0,1,2,\dots,\frac{N-1}{2},-\frac{N-1}{2},\ldots,-2,-1\right],$$ (for $N$ odd), where $l=b-a$ is the extent of the spatial domain; that is $l/N=\Delta x$. (For more details on properties of the DFT and its application to spectral methods see \cite{BriggsHensonBook} and \cite{TrefethenSpectralBook}.)
\begin{equation}
\bcH_{sp}=\sum_j\left[\frac12\pi_j^2+\frac12\left[\mathcal{F}^{-1}_ND_N\mathcal{F}_N\boldsymbol{\varphi}\right]^2_j+\alpha(1-\cos\varphi_j)\right],
\end{equation}
\begin{equation}
\bcS = 
\left(
   \begin{array}{cc}
     0  & \mbox{id} \\
    -\mbox{id} & 0
   \end{array}
\right).
\end{equation}

The resulting system of ODEs is then given by

\begin{equation}
\left[
\begin{array}{c}
{\bf\dot{\boldsymbol{\varphi}}}\\
\bf{\dot{\boldsymbol{\pi}}}
\end{array}
\right]
=\bcS\nabla\bcH_{sp}
=
\left[
\begin{array}{c}
  \boldsymbol{\pi}\\
  -(\mathcal{F}_N^{-1}D_N\mathcal{F}_N)^{\top}(\mathcal{F}_N^{-1}D_N\mathcal{F}_N\boldsymbol{\varphi})
  -\alpha\sin\boldsymbol{\varphi}
\end{array}
\right].
\end{equation}
Again, the integral in the AVF method can be calculated exactly to give
\begin{align}
  \frac{\boldsymbol{\varphi}_{n+1}-\boldsymbol{\varphi}_n}{\Delta t}
  &= (\boldsymbol{\pi}_{n+1}+\boldsymbol{\pi}_n)/2,\\
  \frac{\boldsymbol{\pi}_{n+1}-\boldsymbol{\pi}_n}{\Delta t}
  &=   -(\mathcal{F}_N^{-1}D_N\mathcal{F}_N)^{\top}(\mathcal{F}_N^{-1}D_N\mathcal{F}_N)(\boldsymbol{\varphi}_{n+1}+\boldsymbol{\varphi}_n)/2 \notag\\
&\qquad\qquad\qquad-\alpha(\cos\boldsymbol{\varphi}_{n+1}-\cos\boldsymbol{\varphi}_n)/(\boldsymbol{\varphi}_{n+1}-\boldsymbol{\varphi}_n).
\end{align}

\underline{Initial conditions and numerical data for both discretizations:}

Spatial domain, number $N$ of spatial intervals, and time-step size $\Delta t$ used were\footnote{Here and below, if $x \in [a,b]$, then $\Delta x = \frac{b-a}{N}$, and $x_j = a+j \Delta x, j=0,1,\dots,N$.}

\begin{equation}\notag
x \in [-20,20], \qquad N=200, \qquad \Delta t = 0.01, \qquad \mbox{parameter: } \alpha = 1.
\end{equation}

\begin{figure}[ht]\centering
\includegraphics[scale = 0.3]{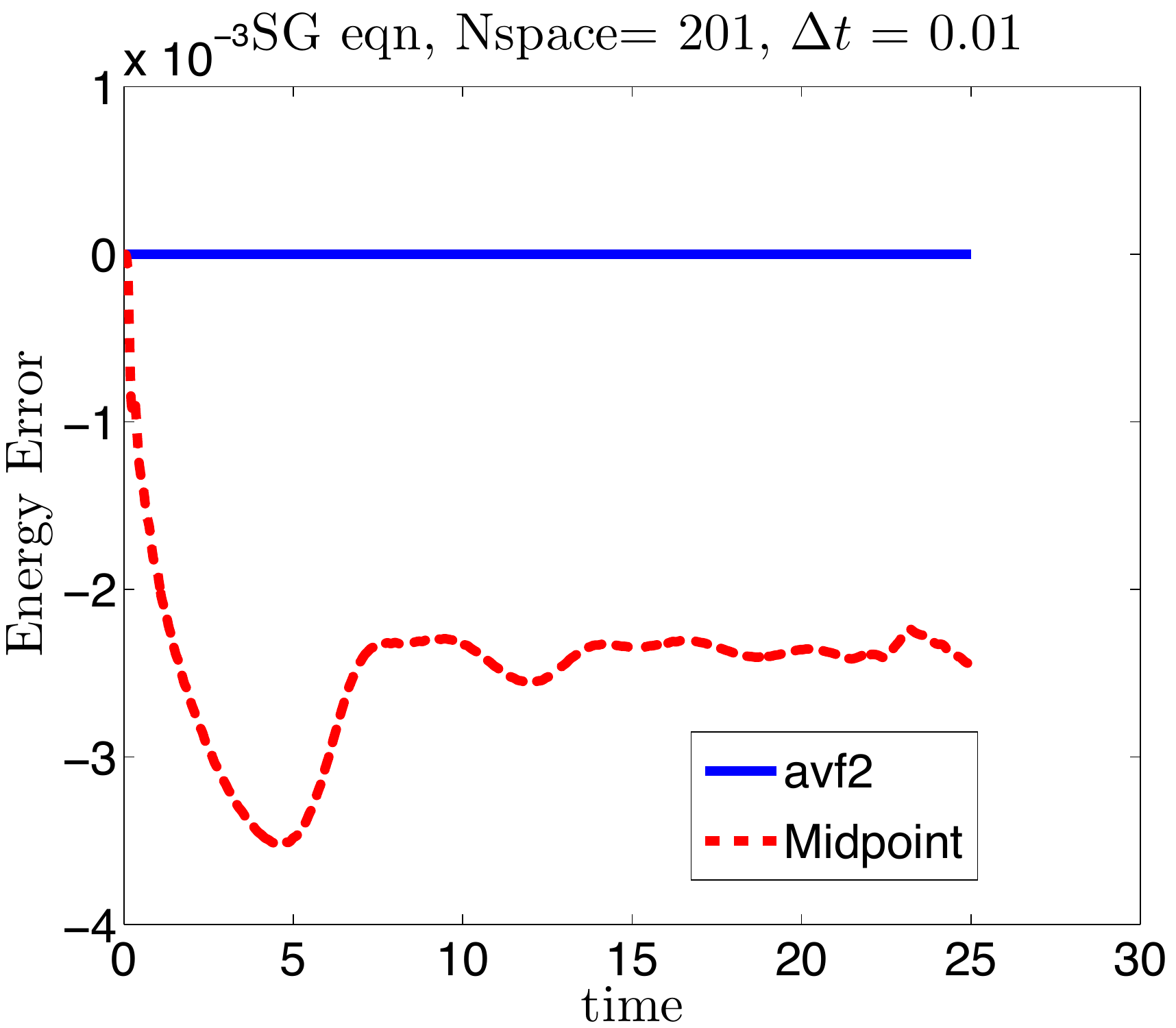}
\hfill
\includegraphics[scale = 0.3]{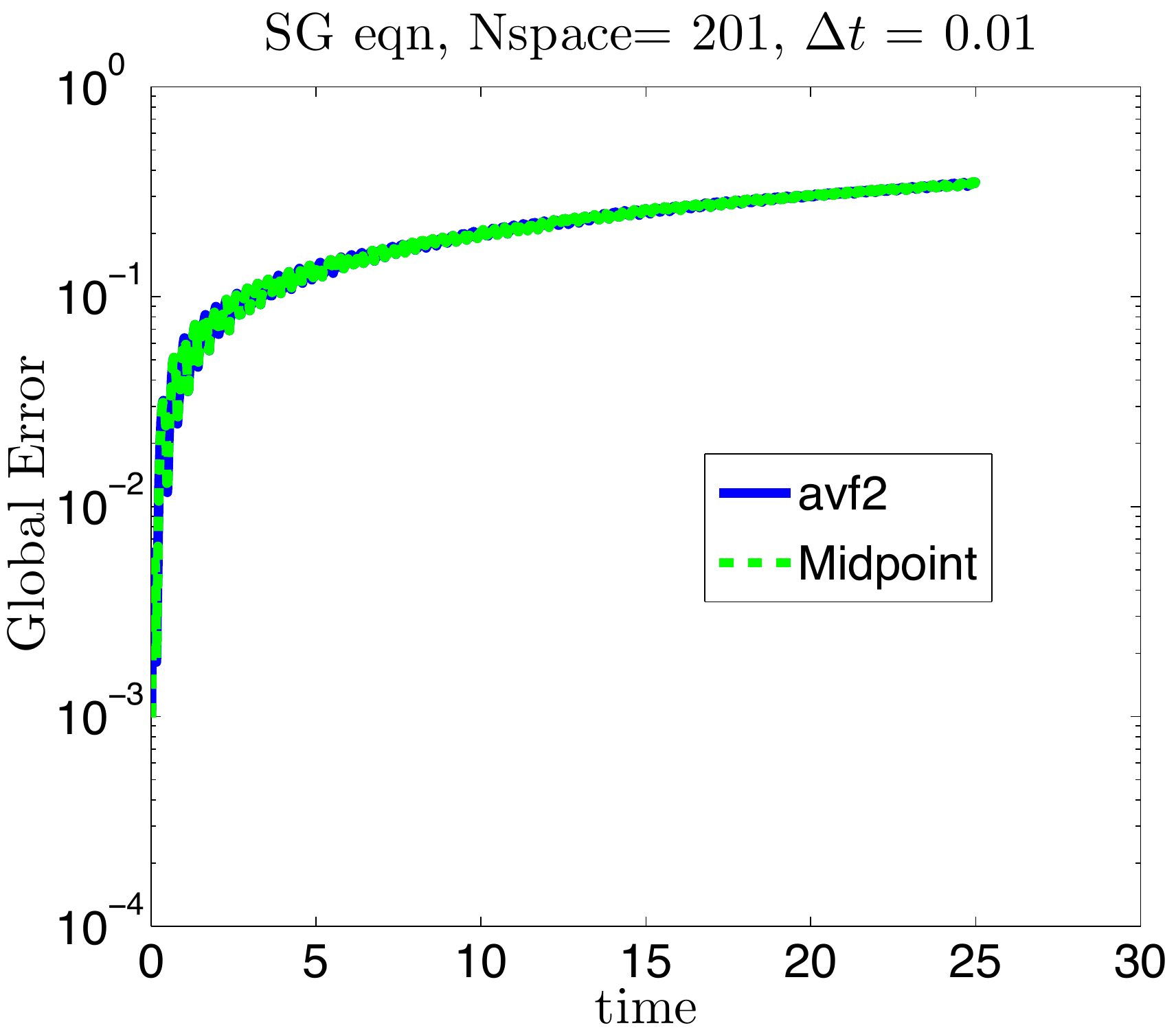}
\caption{Sine-Gordon equation with finite differences semi-discretization: Energy error (left) and global error (right) vs time, for AVF and implicit midpoint integrators.}
\label{Sine-Gordonplot_fd}
\end{figure}

\begin{figure}[ht]\centering
\includegraphics[scale = 0.3]{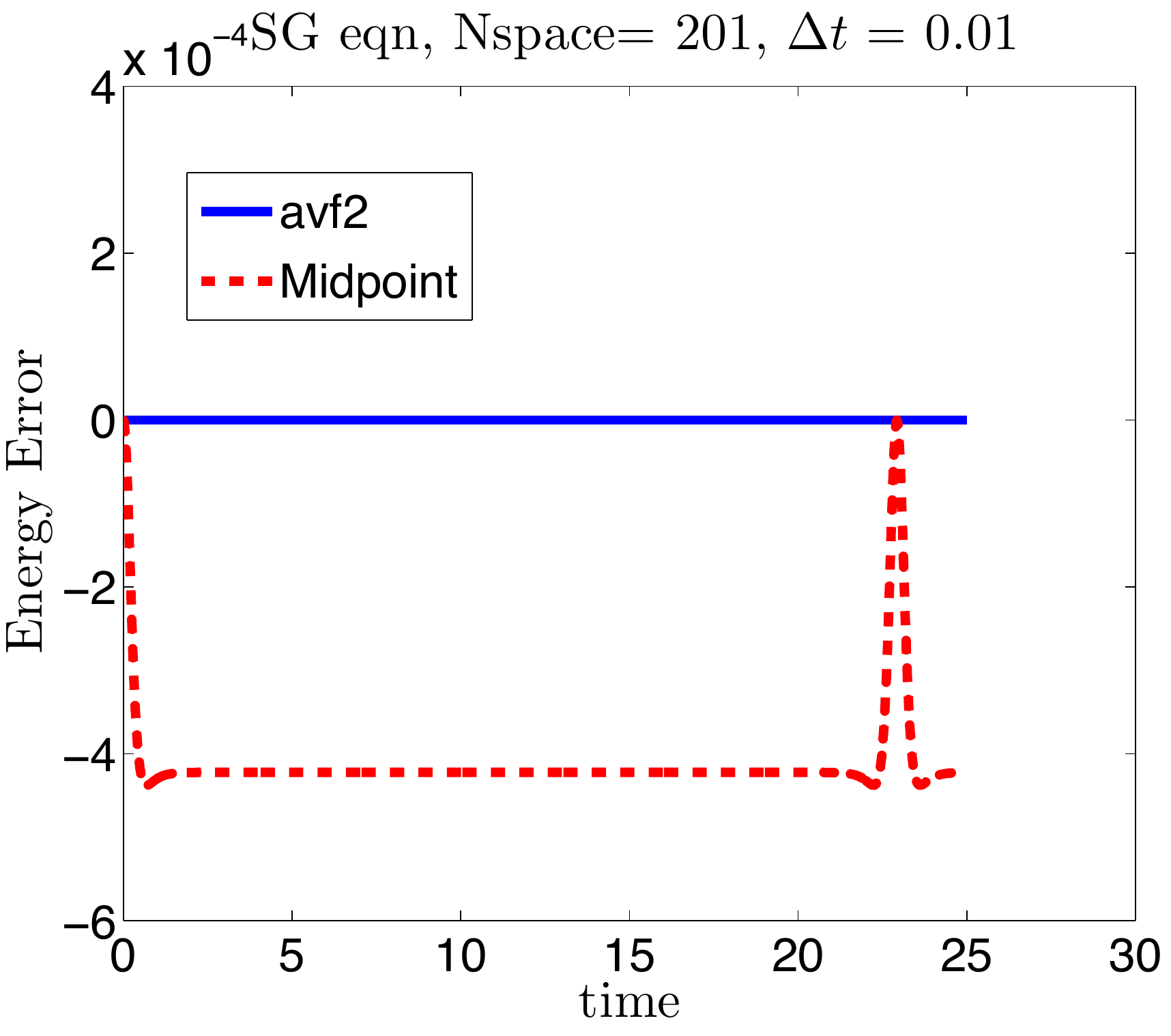}
\hfill
\includegraphics[scale = 0.3]{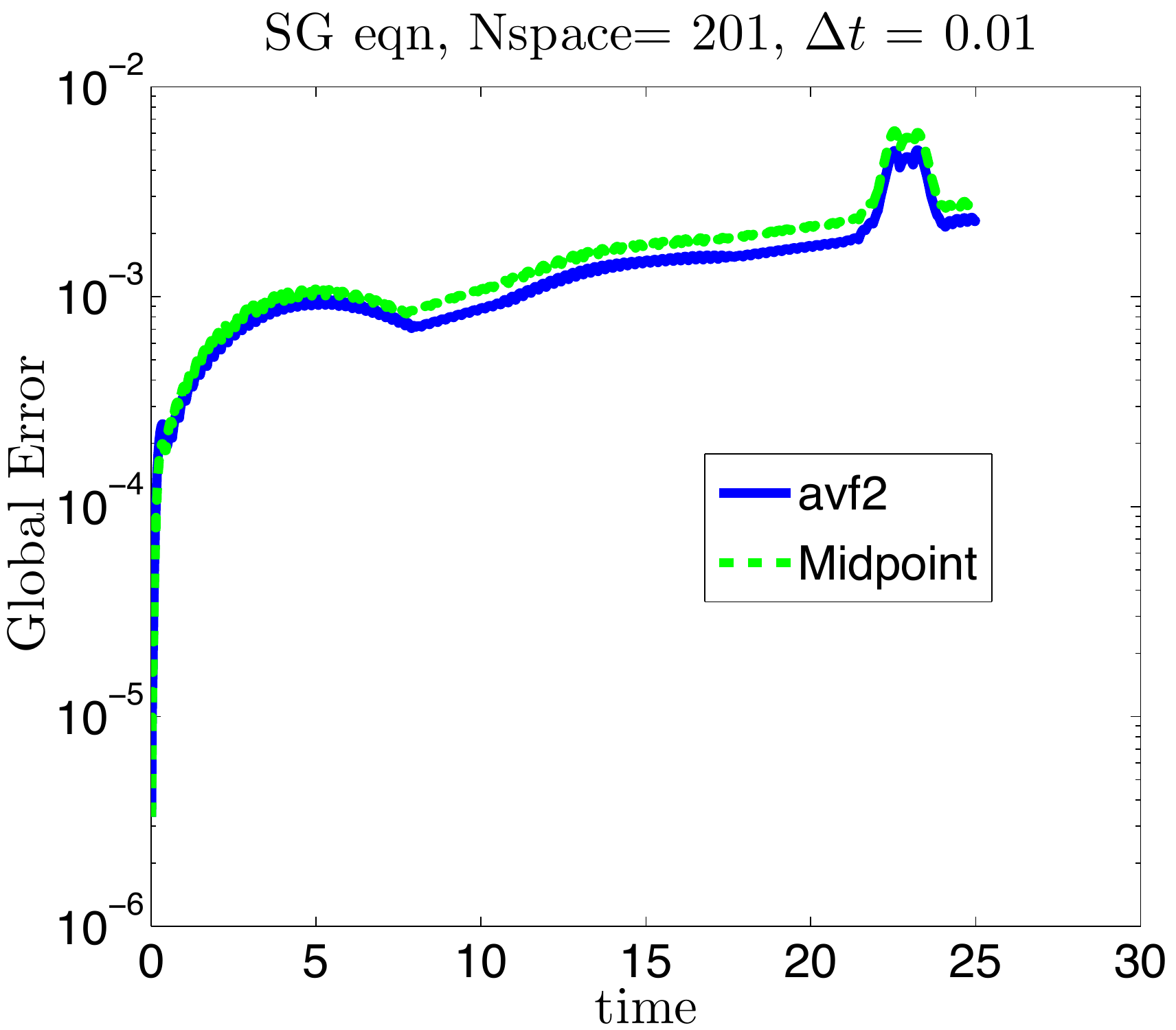}
\caption{Sine-Gordon equation with spectral semi-discretization: Energy error (left) and global error (right) vs time, for AVF and implicit midpoint integrators.}
\label{Sine-Gordonplot_sp}
\end{figure}

Initial conditions:
\begin{equation}
\left.
  \begin{array}{rcl}
     \varphi(x,0)  & = & 0,\\
     \pi(x,0) & = & \dfrac{8}{\cosh (2x)}.  
  \end{array}
\right\}
\qquad \qquad 
\mbox{
   \parbox{3cm}{
       Right-moving kink and left-moving anti-kink solution.
   }
} 
\end{equation}

Numerical comparisons of the AVF method with the well known (symplectic) implicit midpoint integrator\footnote{Recall that the implicit midpoint integrator is given by $\frac{u_{n+1}-{u_n}}{\Delta t}=f\left(\frac{u_n+u_{n+1}}{2}\right)$.} are given in figure \ref{Sine-Gordonplot_fd} for the finite differences discretization, and in figure \ref{Sine-Gordonplot_sp} for the spectral discretization.

\begin{example} 
Korteweg--de Vries equation:
\end{example}

\underline{Continuous:}

\begin{equation}
\frac{\partial u}{\partial t} = -6u\frac{\partial u}{\partial x} - 
\frac{\partial^3 u}{\partial x^3},
\end{equation}

\begin{equation}
\cH = \int 
  \left[
     \frac{1}{2} \left( u_x \right)^2 - u^3
  \right] \, dx,
\end{equation}
\begin{equation}
\cS = \frac{\partial}{\partial x}.
\end{equation}
Boundary conditions: periodic, $u(-20,t) = u(20,t)$.

\underline{Semi-discrete: }
\begin{equation}
\bcH =  \sum_j 
  \left[
     \frac{1}{2 (\Delta x)^2}\left( u_{j+1}-u_{j} \right)^2 - u_j^3
  \right],
\end{equation}
\begin{equation}
\bcS = \frac{1}{2\Delta x}
\left[
   \begin{array}{ccccc}
   0 &   -1   &        &        &  1 \\
   1 &    0   &   -1   &        &    \\
     & \ddots & \ddots & \ddots &    \\
     &        &    1   &    0   & -1 \\ 
  -1 &        &        &    1   &  0 
   \end{array}
\right].
\end{equation}

\underline{Initial conditions and numerical data:}

\begin{equation}\notag
x \in [-20,20], \qquad N=400, \qquad  \Delta t
=0.001.
\end{equation}

\begin{figure}\centering
\includegraphics[scale = 0.45]{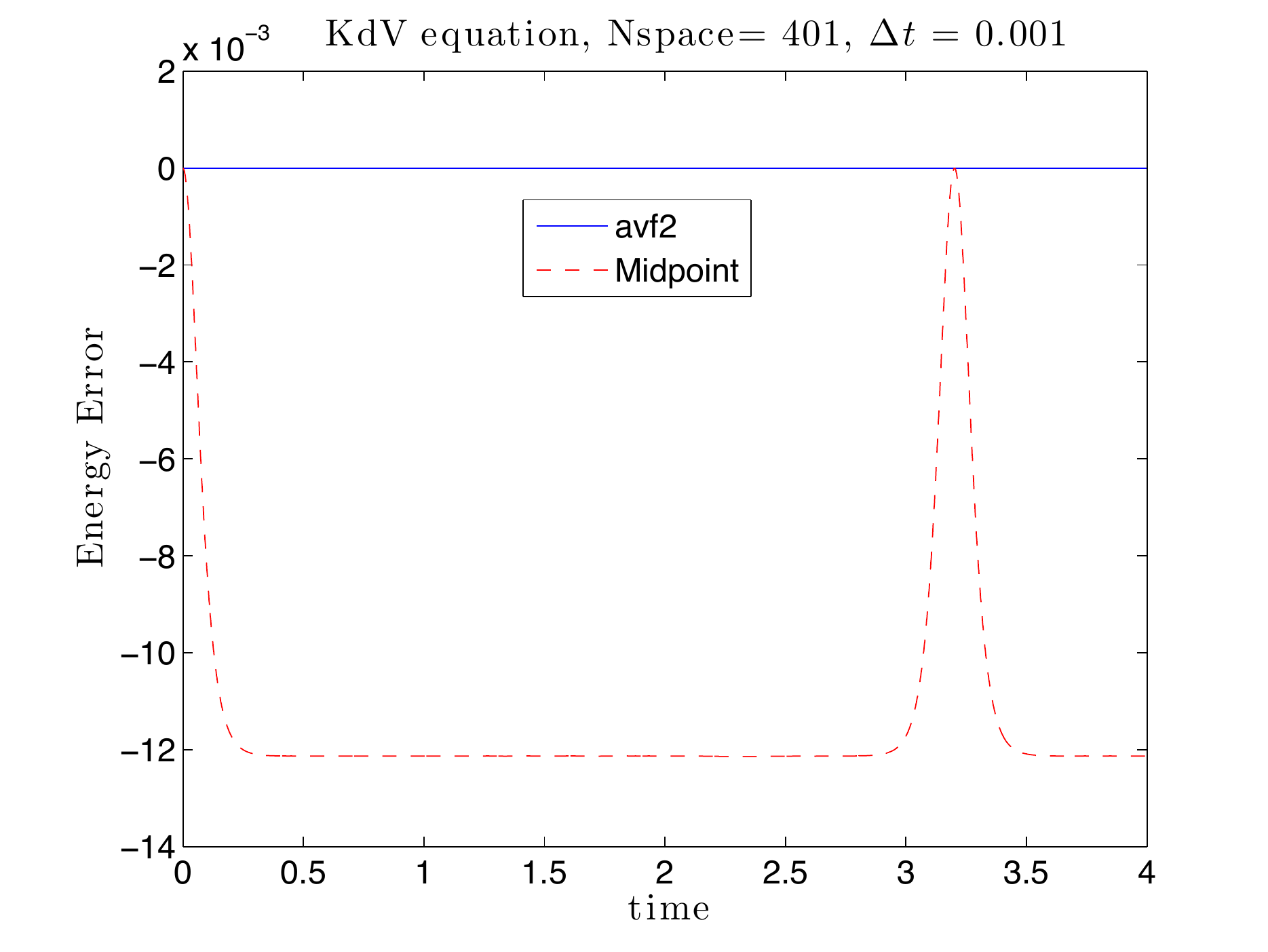}
\includegraphics[scale = 0.45]{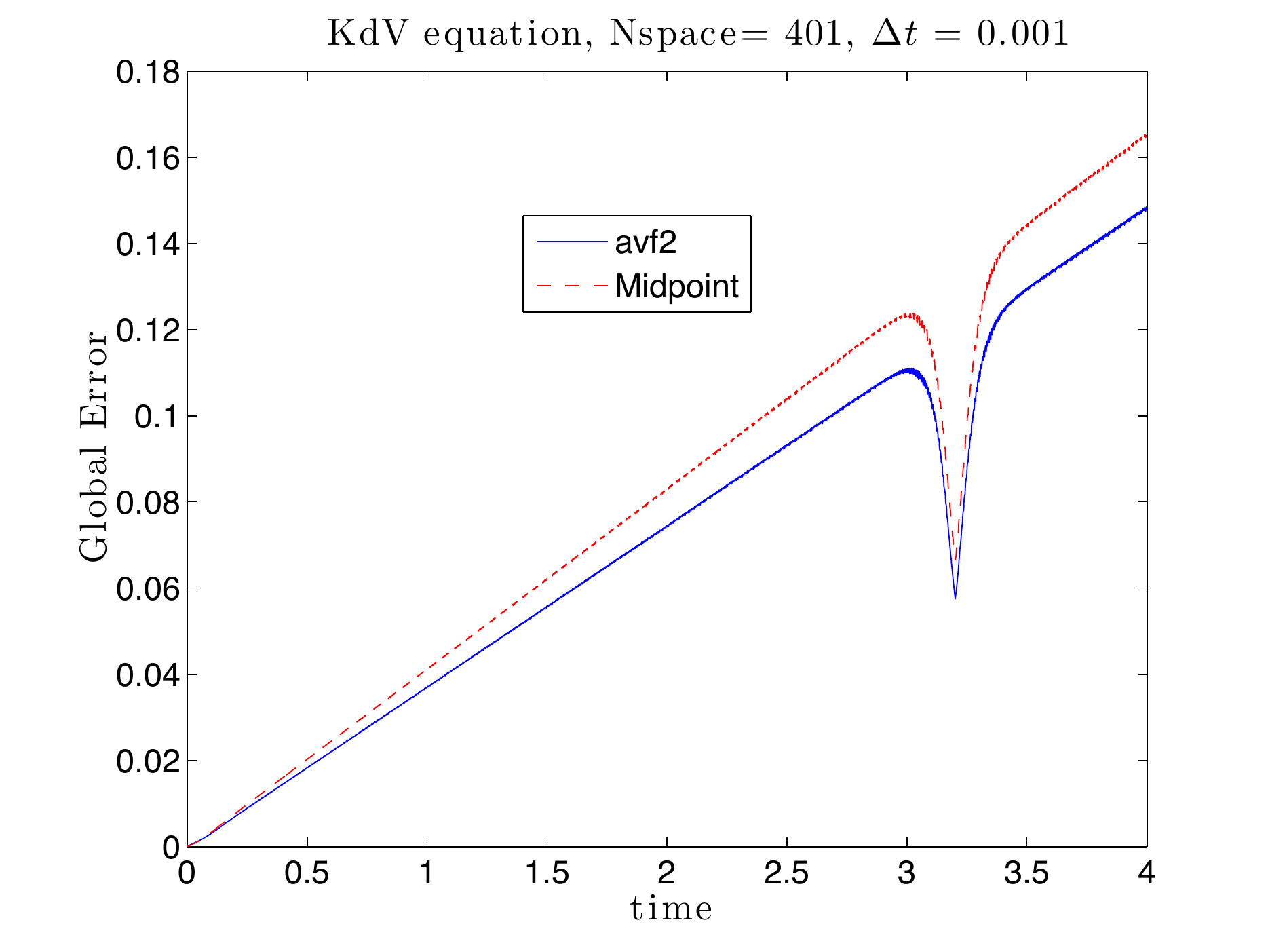}
\caption{Korteweg--de Vries equation: Energy error (left) and global error (right) vs time, for AVF and implicit midpoint integrators.}
\label{KdVplot}
\end{figure}

Initial condition:  $u(x,0)  =  6 \; \mbox{sech}^2(x)$ (for two solitons). Numerical comparisons of the AVF and the midpoint rule are give in Figure \ref{KdVplot}. The global errors are comparable. The AVF is particularly easy to implement here because the integral (\ref{avf}) is evaluated exactly by Simpson's rule.
\bigskip
\begin{example} 
Nonlinear Schr\"{o}dinger equation:
\end{example}

\underline{Continuous:}

\begin{equation}
  \frac{\partial}{\partial t}
  \left(
     \begin{array}{c}
       u \\ u^*
     \end{array}
  \right)
  =
  \left(
     \begin{array}{cc}
       0 & i \\
      -i & 0
     \end{array}
  \right)
  \left(
     \begin{array}{c}
       \frac{\delta \cH}{\delta u} \\ \frac{\delta \cH}{ \delta u^*}
     \end{array}
  \right),
\end{equation}
where $u^*$ denotes the complex conjugate of $u$.

\begin{equation}
  \cH = \int 
  \left[
    -\left| \frac{\partial u}{\partial x} \right|^2 
    + \frac{ \gamma}{2} |u|^4
  \right] \, dx,
\end{equation}
\begin{equation}
  \cS=
  \left(
     \begin{array}{cc}
       0 & i \\
      -i & 0
     \end{array}
  \right).
\end{equation}

Boundary conditions: periodic, $u(-20,t) = u(20,t)$.

\underline{Semi-discrete: }
\begin{equation}
  \bcH = \sum_j 
  \left[
     -\frac{1}{(\Delta x)^2} \left| u_{j+1}-u_{j} \right|^2 + \frac{\gamma}{2} |u_j|^4
  \right],
\end{equation}
\begin{equation}
  \bcS = 
 i  \left(
     \begin{array}{cc}
             0     & \mbox{id} \\
        -\mbox{id} &     0 
     \end{array}
  \right).
\end{equation}

\underline{Initial conditions and numerical data:}

\begin{equation} \notag
x \in [-20,20], \qquad N=200,  \qquad \Delta t = 0.05, \qquad \mbox{parameter:} \, \gamma=1.
\end{equation}

\begin{figure}\centering
\includegraphics[scale = 0.45]{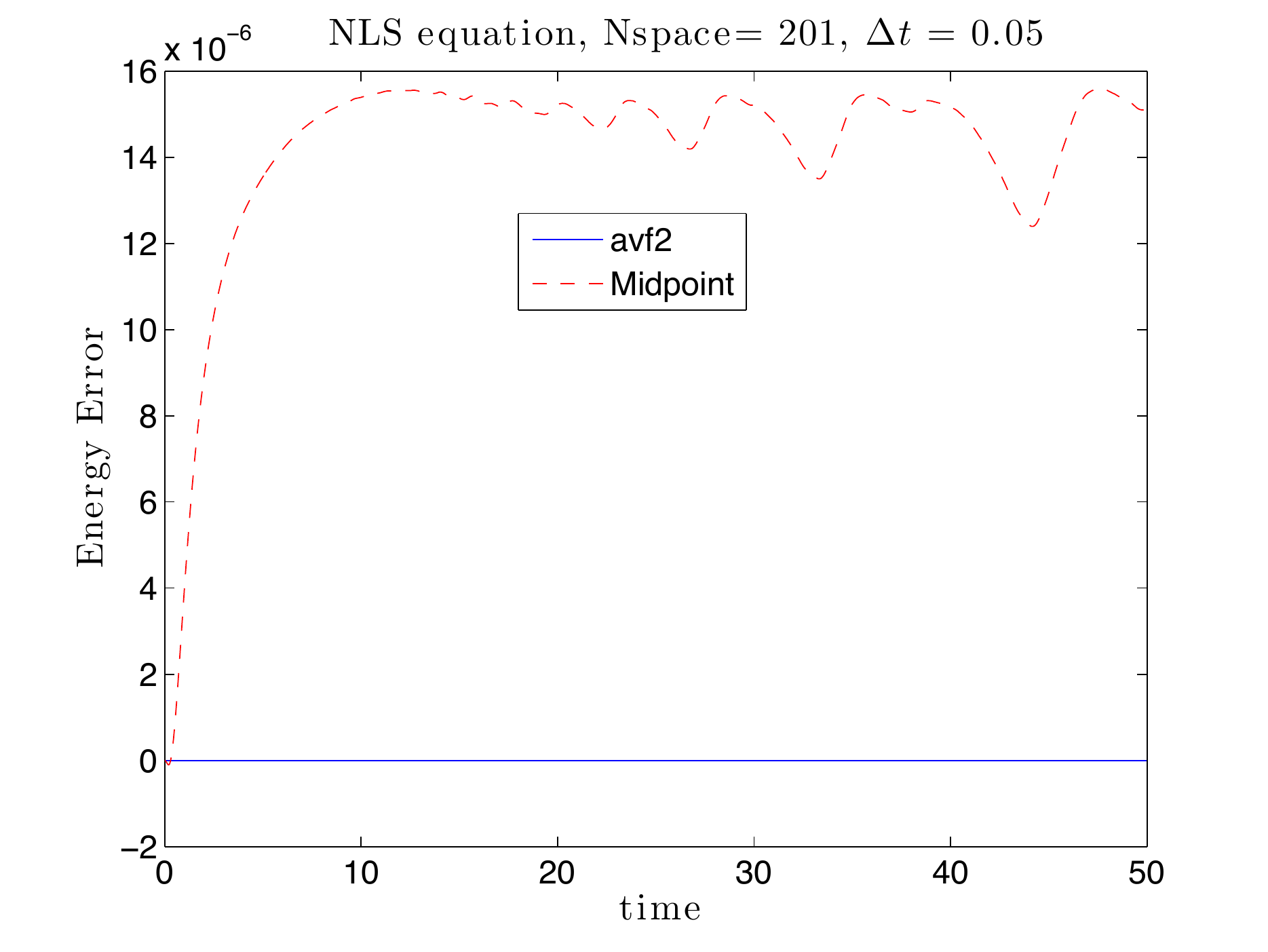}
\hfill
\includegraphics[scale = 0.45]{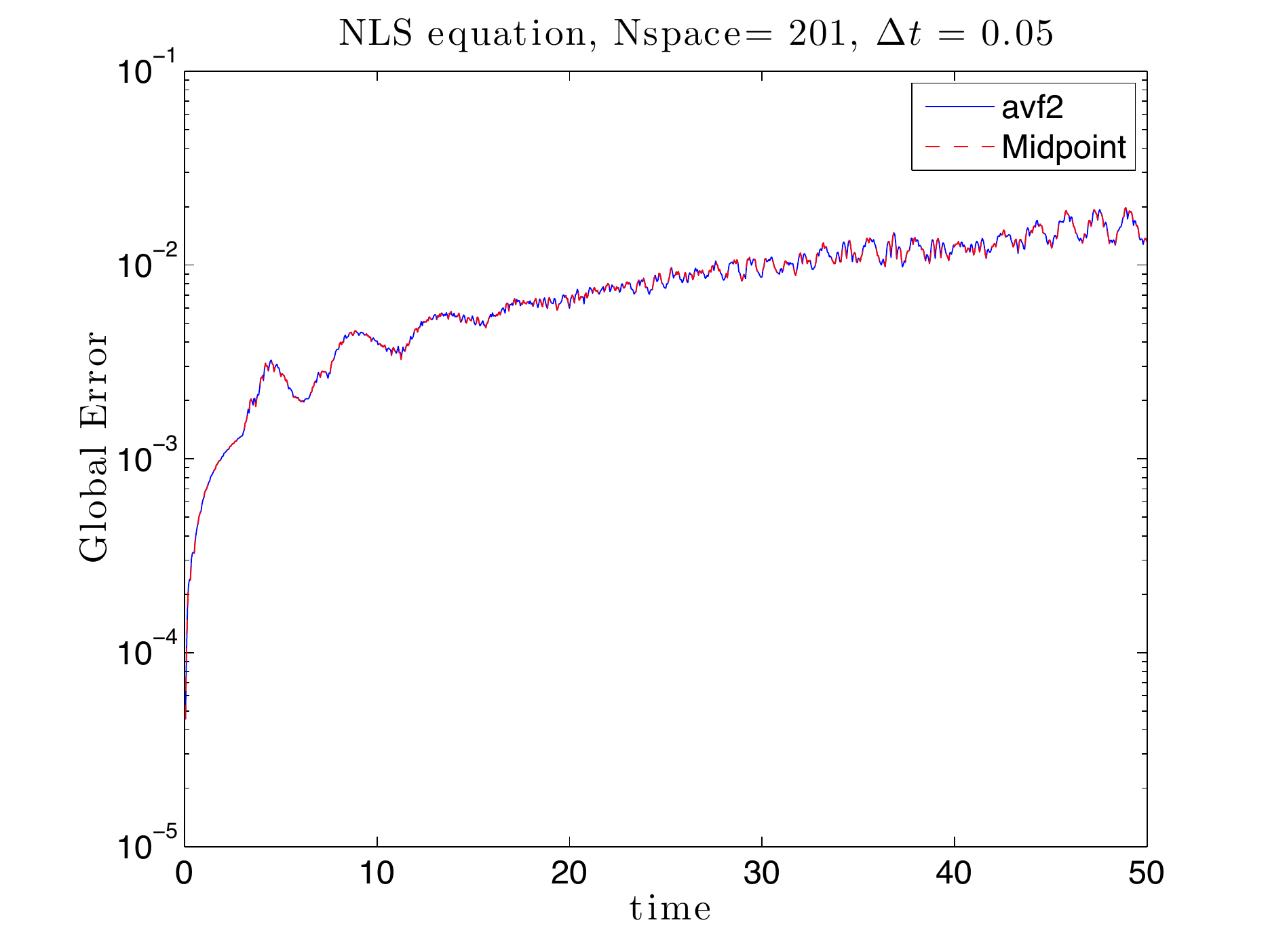}
\caption{Nonlinear Schr\"{o}dinger equation: Energy  and global error vs time, for AVF and implicit midpoint integrators.}
\label{NLSplot1}
\end{figure}

\begin{figure}\centering
\includegraphics[scale = 0.45]{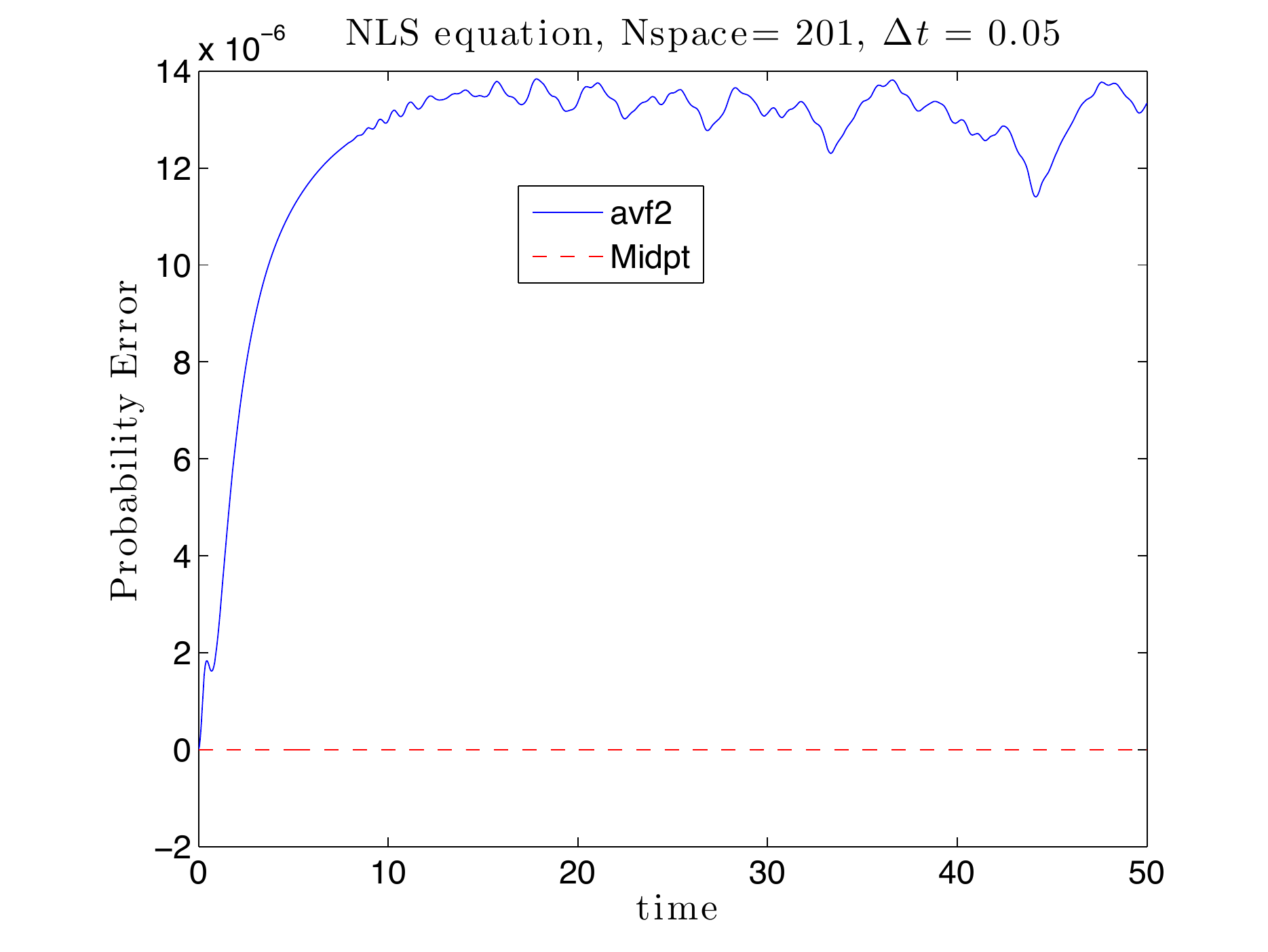}
\caption{Nonlinear Schr\"{o}dinger equation: Total probability error vs time, for AVF and implicit midpoint integrators.}
\label{NLSplot2}
\end{figure}

Initial conditions:

\begin{equation}
  \left\{
    \begin{array}{rcl}
      \Re  u(x,0) &=& \mbox{exp} \left( {-(x-1)^2\slash 2} \right), \\
      \Im  u(x,0) &=& \mbox{exp} \left( {- x^2\slash 2} \right).
    \end{array}
  \right.
\end{equation}
Numerical comparisons of the AVF and the midpoint rule are give in Figures \ref{NLSplot1}, \ref{NLSplot2}. The global errors are comparable. As for the KdV equation,  the integral (\ref{avf}) is evaluated exactly by Simpson's rule. The AVF method does not conserve the total probalility because it is not unitary.

\begin{example}
Nonlinear Wave Equation:
\end{example}

\underline{Continuous:}

The 2D wave equation
\begin{equation}
\label{2dwavepde}
\frac{\partial^2 \varphi}{\partial t^2} = \Delta \varphi -
\frac{\partial V(\varphi)}{\partial \varphi},\quad \varphi=\varphi(x,y,t),\quad (x,y)\in [-1,1]\times [-1,1], \, t\ge0,
\end{equation}
is a Hamiltonian PDE with Hamiltonian function
\begin{equation}
\mathcal{H} = \int_{-1}^{1}\int_{-1}^{1}
\left[ 
   \frac{1}{2}(\pi^2+\varphi_x^2+\varphi_y^2)+V(\varphi)
\right] \,dx\,dy ,
\end{equation}
where $\pi=\partial / \partial t \varphi$ and the operator $\mathcal{S}$ is the canonical $2\times 2$ symplectic
matrix.

Boundary conditions: periodic.

\underline{Semi-discrete: spectral elements}

We
discretize the Hamiltonian in space with a tensor product Lagrange
quadrature formula based on
$p+1$ Gauss-Lobatto-Legendre (GLL) quadrature nodes in each space direction.
We obtain
\begin{equation}
\bcH = \frac{1}{2}\sum_{j_1=0}^p\sum_{j_2=0}^p w_{j_1}
w_{j_2}\!\left(\!\pi_{j_1,j_2}^2 \! +\!
\left(\sum_{k=0}^pd_{j_1,k}\varphi_{k,j_2}\right)^2\!\!+\!\left(\sum_{m=0}^pd_{j_2,m}\varphi_{j_1,m}\!\right)^2\!\! +\! \frac{1}{2}\varphi_{j_1,j_2}^4\right),
\end{equation}
where $d_{j_1,k}=\left.\frac{d l_{k}(x)}{dx}\right|_{x=x_{j_1}}$, and
$l_{k}(x)$ is the $k$-th Lagrange
basis function based on the GLL quadrature nodes $x_0,\dots , x_p$, and with
$w_0,\dots , w_p$ 
the corresponding quadrature weights. The numerical approximation is
\begin{equation}
\varphi_p(x,y,t)=\sum_{k=0}^p\sum_{m=0}^p \varphi_{k,m}(t)l_k(x)l_m(y),
\end{equation}
and has the property $\varphi_p(x_{j_1},y_{j_2},t)=\varphi_{j_1,j_2}(t)$, so
that the data can be stored in the
$(p+1)\times (p+1)$  matrix with entries $\varphi_{j_1,j_2}$.

\underline{Initial conditions and numerical data:}
$$
(x,y) \in [-1,1]^2, \quad V(\varphi) = \frac{\varphi^4}{4}.
$$
Initial condition: $\varphi(x,y,0) = \mathrm{sech}(10 x)\mathrm{sech}(10
y)$, $\pi (x,y,0)=0$. 

\begin{figure}
\begin{center}
\begin{tabular}{cc}
$T=0$ & $T=3.1250$ \\
\includegraphics[width=4.5cm]{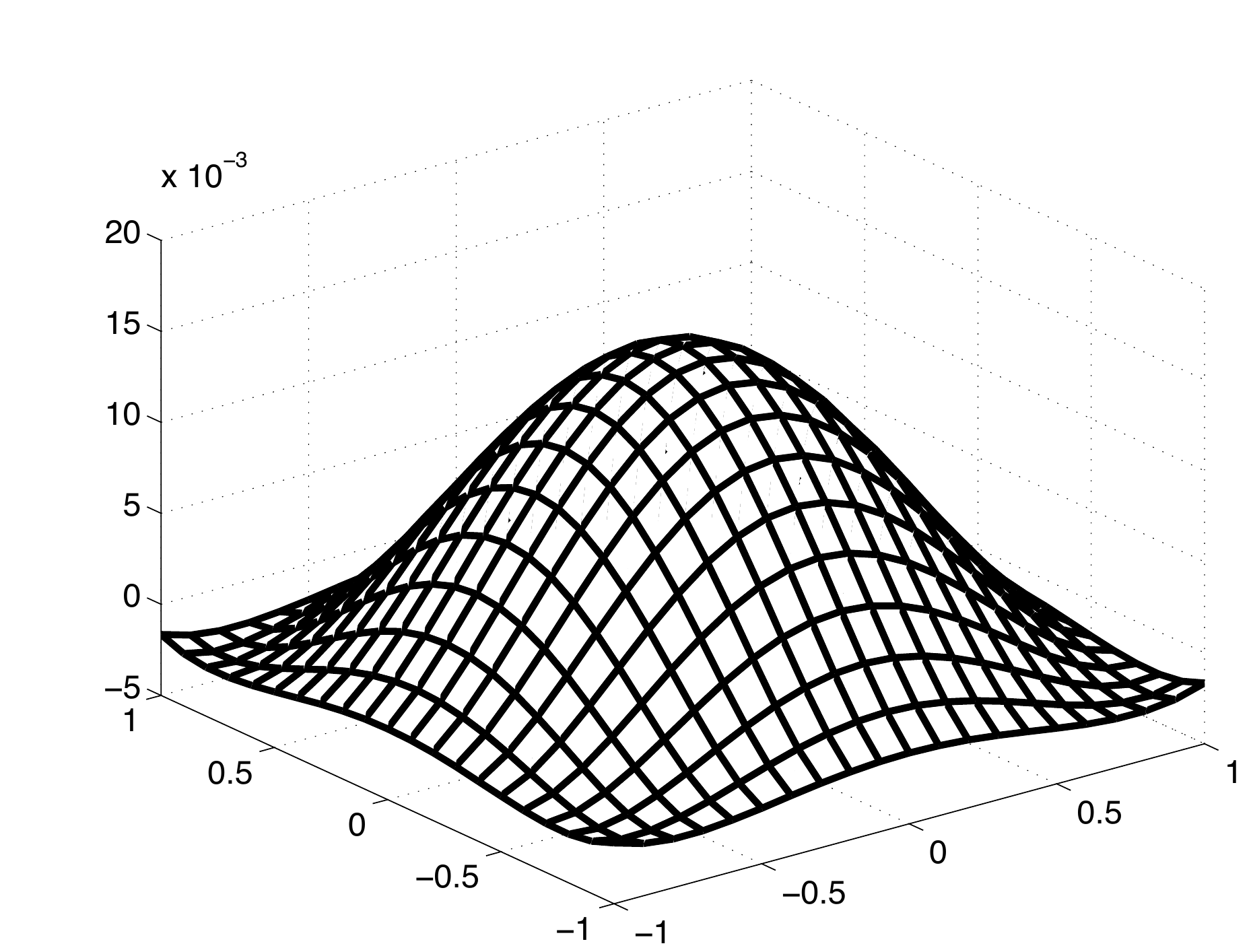} &
\includegraphics[width=4.5cm]{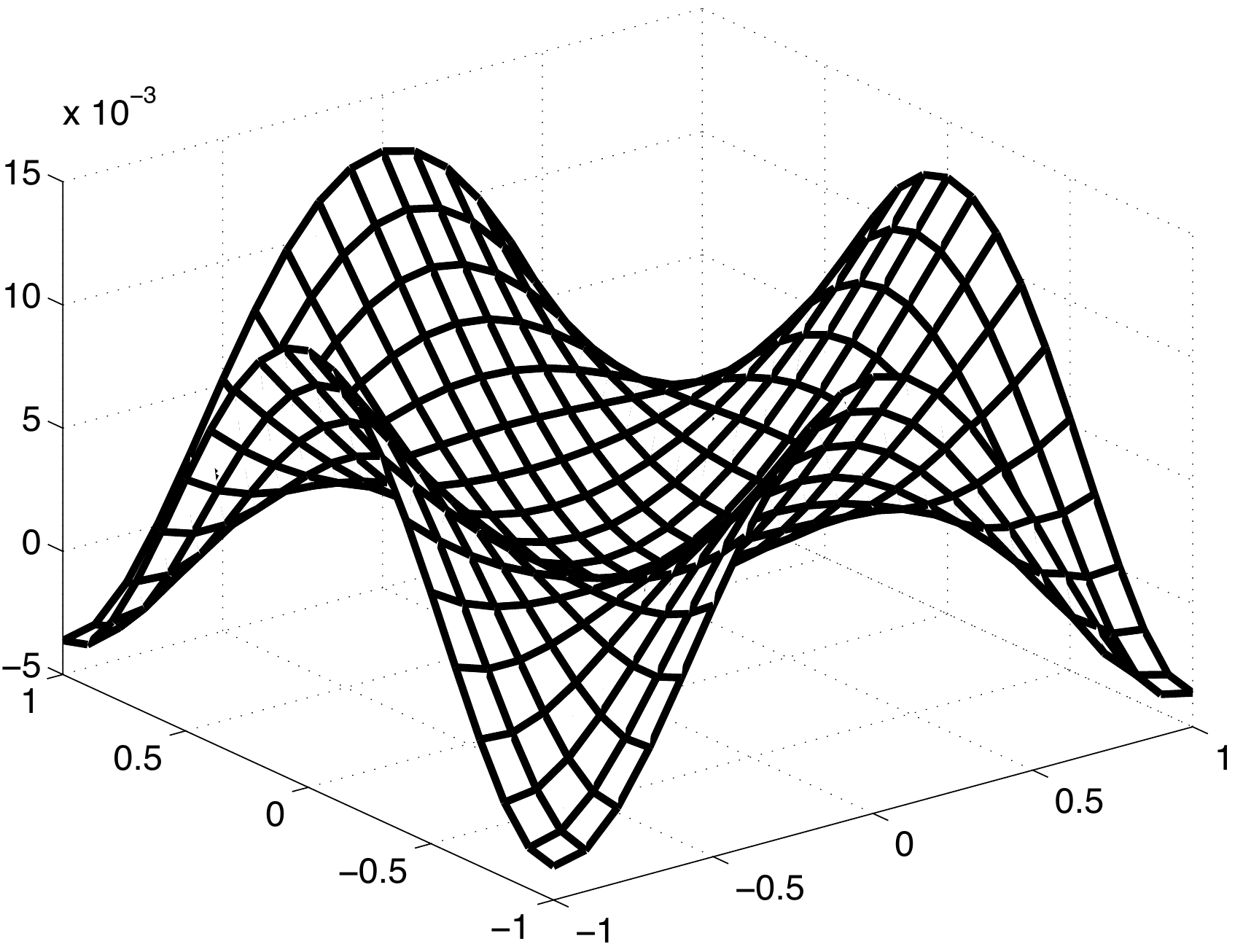}\\
$T=6.875$ & $T=10$ \\
\includegraphics[width=4.5cm]{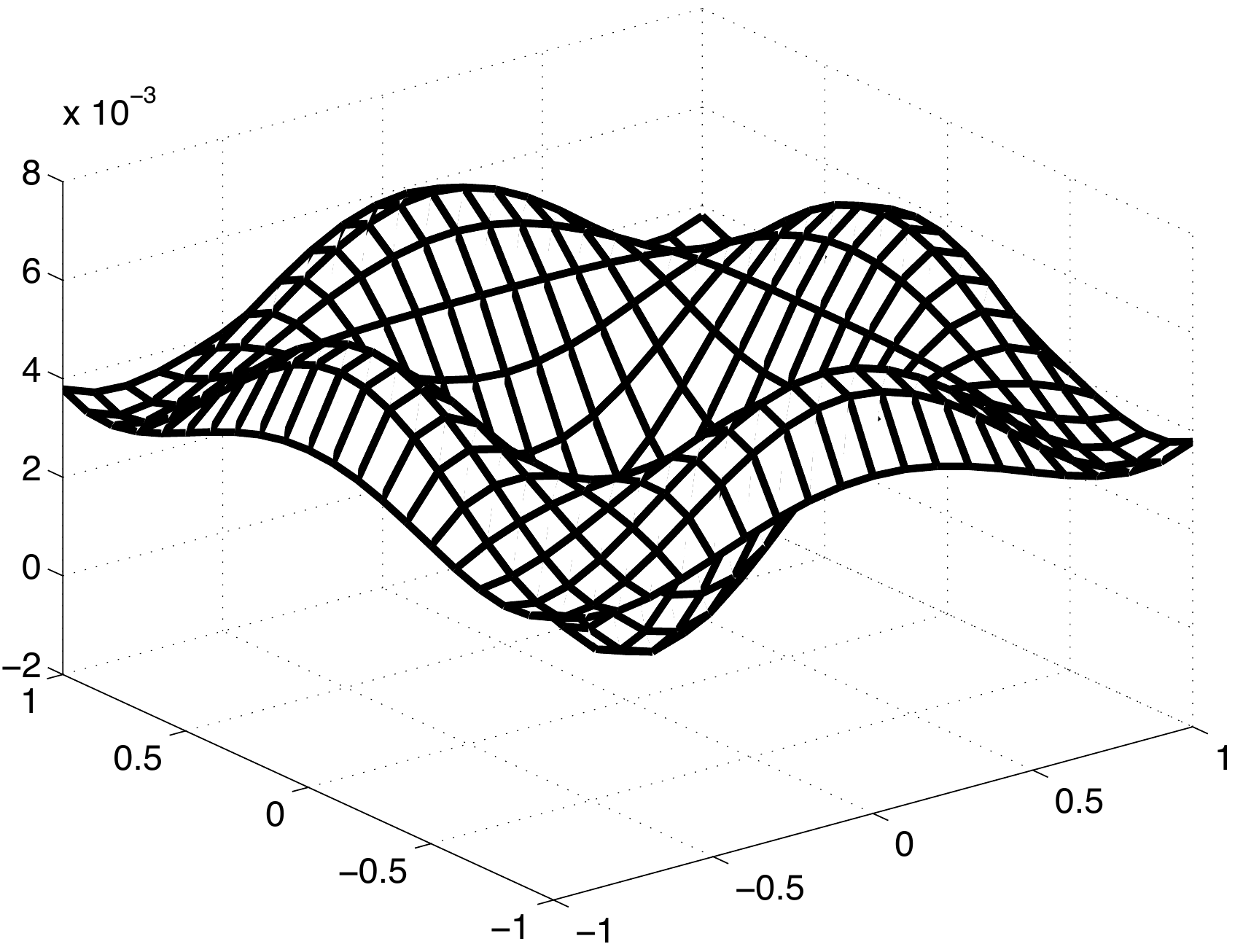} &
\includegraphics[width=4.5cm]{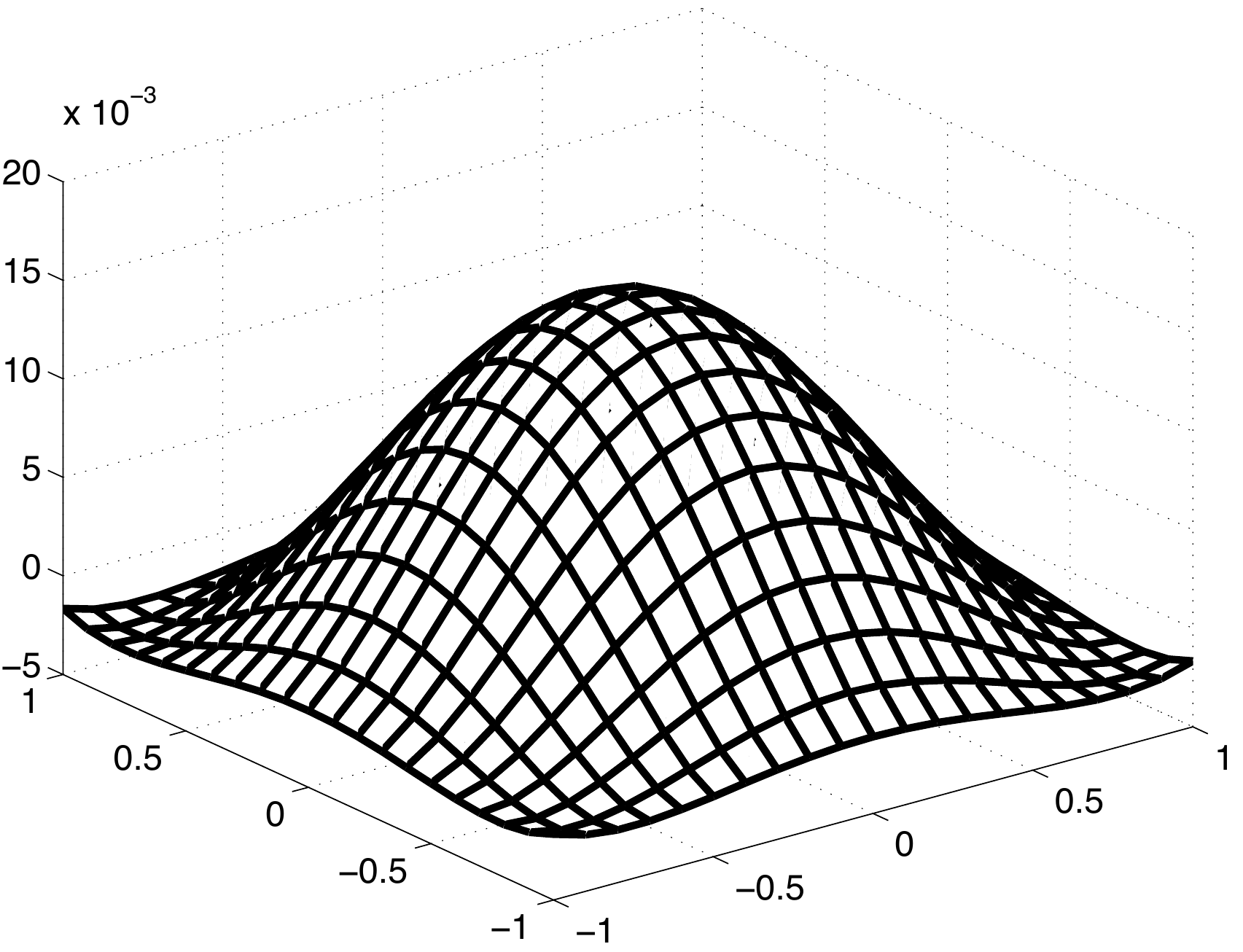}\\
\end{tabular}
\caption{Snapshots of the solution of the 2D wave equation at different
times. AVF method with step-size $\Delta t=0.6250$. Space
discretization with $6$ Gauss Lobatto nodes in each space direction. Numerical solution interpolated on a equidistant grid of $21$ nodes in each space direction. \label{fig:wave}}
\end{center}
\end{figure}

\begin{figure}
\begin{center}
\includegraphics[scale = 0.45]{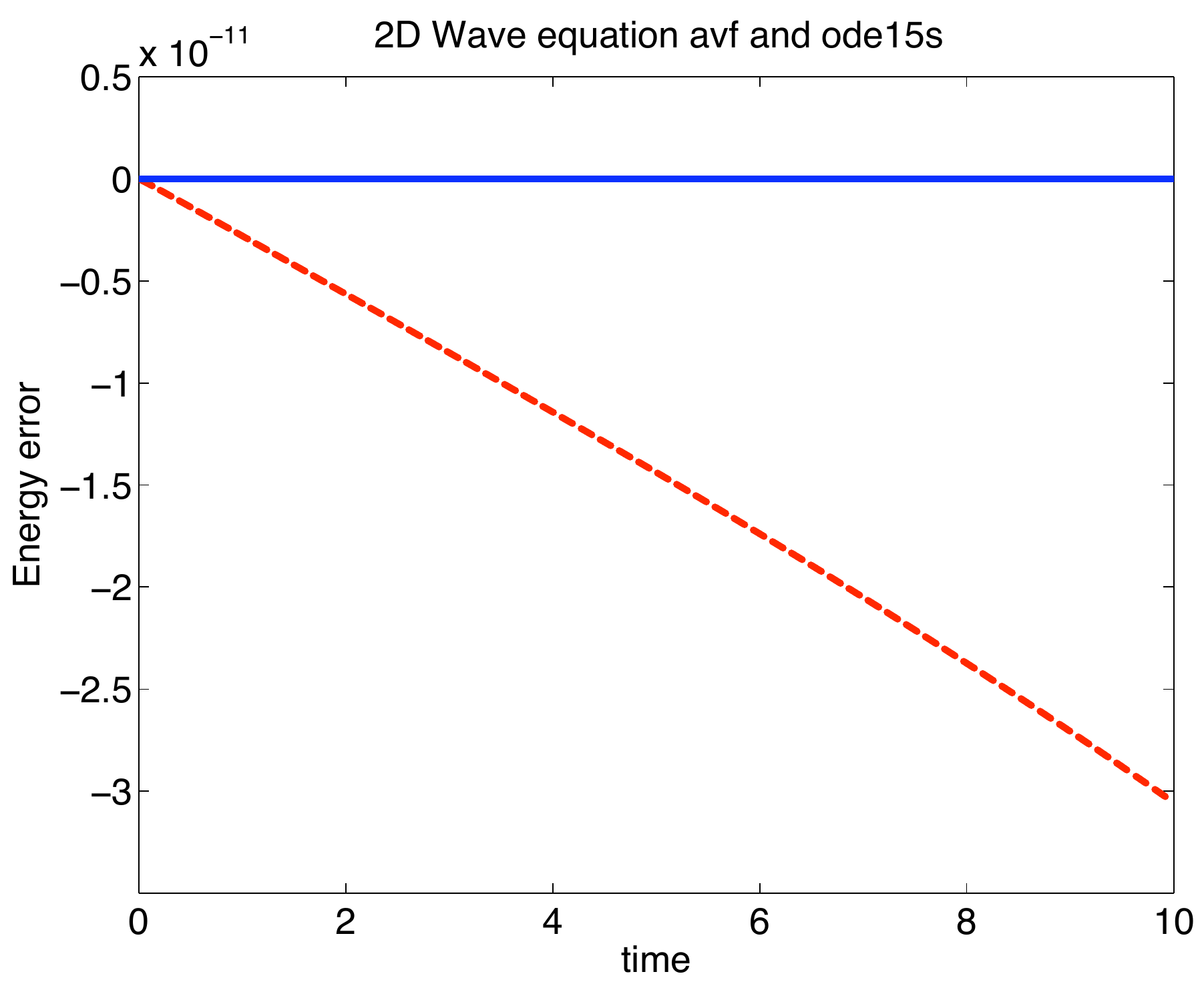} \\
\caption{The 2D wave equation (\ref{2dwavepde}). MATLAB routine \texttt{ode15s} with absolute and relative tolerance
$10^{-14}$ (dashed line), and AVF method with step size $\Delta t=10/(2^5)$ (solid
line). Energy error versus time. Time interval $[0,10]$. Space
discretization with $6$ Gauss Lobatto nodes in each space direction.
\label{fig:energywave}}
\end{center}
\end{figure}

\medskip
In figure~\ref{fig:wave} we show some snapshots of the solution. The energy
error is shown in figure~\ref{fig:energywave}.

\subsection{Linear conservative PDEs}
\begin{example}
(Linear) Time-dependent Schr\"{o}dinger Equation:
\end{example}

\underline{Continuous:}

\renewcommand{\arraystretch}{1.2}
\begin{equation}
  \frac{\partial u}{\partial t}
  = i \frac{\partial^2 u}{\partial x^2} - i V(x) u.
\end{equation}
\renewcommand{\arraystretch}{1}

This equation is bi-Hamiltonian, i.e. it has 2 independent symplectic structures. The first Hamiltonian formulation has
\begin{equation}
  \cH_1 = \int_{-\pi}^{\pi} 
  \left[
    - \left| \frac{\partial u}{\partial x} \right|^2 - V(x) \left| u \right|^2 
  \right] \, dx
\end{equation}
and
\begin{equation}
  \cS_1 =
  \left(
     \begin{array}{cc}
       0 & i \\
       -i & 0
     \end{array}
  \right).
\end{equation}
The second Hamiltonian formulation has
\begin{equation}
  \cH_2 = \int_{-\pi}^{\pi}  |u|^2 \, dx
\end{equation}
and
\begin{equation}
  \cS_2 =
  \left(
     \begin{array}{cc}
       0 & \partial_x^2 - V(x) \\
      - \partial_x^2 + V(x) & 0
     \end{array}
  \right).
\end{equation}

Boundary conditions: periodic, $u(-\pi,t) = u(\pi,t)$.

\underline{Semi-discrete: }
\begin{equation}
  \bcH_1 = \sum_{j}^{} 
  \left[
    -\frac{1}{(\Delta x)^2} \left| u_{j+1}-u_{j} \right|^2 - V(x_j)  |u_j|^2
  \right],
\end{equation}
\begin{equation}
  \bcS_1 = 
 i  \left(
     \begin{array}{cc}
             0     & \mbox{id} \\
        -\mbox{id} &     0 
     \end{array}
  \right).
\end{equation}
The second semi-discretization is
\begin{equation}
 \bcH_2 = \sum_{j}^{} |u_j|^2,
\end{equation}
\begin{equation}
  \bcS_2 = 
 i  \left(
     \begin{array}{cc}
             0     & A  \\
        - A &    0
     \end{array}
  \right),
\end{equation}
where
\begin{equation}
  A = 
  \left(
     \begin{array}{ccccc}
            -2-V & 1      & 0 & \dots & 1 \\
               1  &  -2-V & 1 &          & 0 \\
               0 & \ddots & \ddots & \ddots & \vdots \\
               0 & \dots &  &  & 1  \\
               1 & 0 & \dots & 1 & -2-V
     \end{array}
  \right).
\end{equation}
Both discretizations result in the same semi-discrete system and the AVF method (which in the linear case coincides with the midpoint rule) therefore preserves both $\bcH_1$ and $\bcH_2$, as well as the two symplectic structures.

\medskip 
\underline{Initial conditions and numerical data:}

\begin{equation} \notag
x \in [-\pi,\pi], \qquad N=50,  \qquad \Delta t = 0.1, \qquad 
V(x)=1-\cos(x).
\end{equation}

Initial conditions:
\begin{equation}
\Re u(x,0)=e^{-(\frac{x}{2})^2}, \qquad \Im u(x,0)=0.
\end{equation}

\begin{figure}\centering
\includegraphics[scale = 0.45]{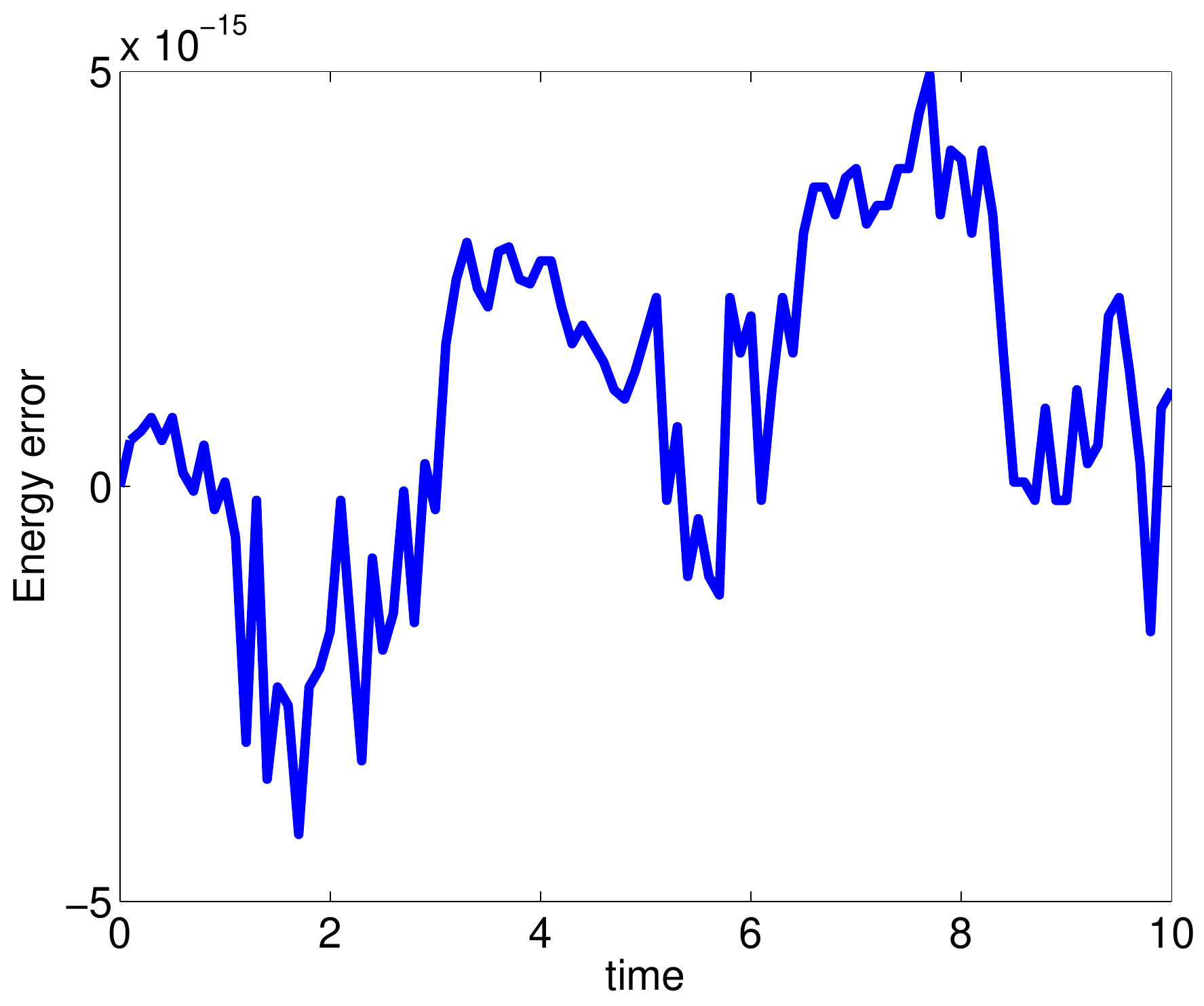}
\caption{Linear Schr\"{o}dinger equation: Error in energy $\bcH_1 \Delta x$ vs time, AVF method}
\label{LSplot}
\end{figure}

\begin{example} 
Maxwell's Equations (1d):
\end{example}
We first look at the one-dimensional Maxwell equation

\underline{Continuous:}

\renewcommand{\arraystretch}{1.2}

\begin{equation}
 \frac{\partial}{\partial t} \left[
    \begin{array}{c}
      E \\
      B
    \end{array}
  \right]  
  =
  \left[
    \begin{array}{cc}
      0 & c \frac{\partial}{\partial x} \\
      c \frac{\partial}{\partial x} & 0 
    \end{array}
  \right]
  \left[
    \begin{array}{c}
     \frac{\delta \cH}{\delta E} \\
     \frac{\delta \cH}{\delta B} 
    \end{array}
  \right],
\end{equation}

\renewcommand{\arraystretch}{1.0}

\begin{equation}
\cH (E,B)= \int_0^1 c\frac{1}{2} \left( E^2 + B^2 \right) \, dx,
\label{integral}
\end{equation}
and 
\begin{equation}
\cS=\left[
    \begin{array}{cc}
      0 & \frac{\partial}{\partial x} \\
      \frac{\partial}{\partial x} & 0 
    \end{array}
  \right]. 
\end{equation}
$\cS$ is skew-adjoint on $\{ (E,B) \in C^1 : E(0)=E(1)=0 \}$ (and
therefore on the Sobolev space $H_0^1$). 

Boundary conditions:
\begin{equation}
 \left\{
  \begin{array}{c}
  E(0,t)=E(1,t)=0, \\
  \frac{\partial B}{\partial x}(0,t) = \frac{\partial B}{\partial x}(1,t) = 0.
  \end{array}
  \right.
\end{equation}

\underline{Semi-discrete: }

We now obtain $\bcH$ 
by discretizing $\cH$ in a simple way by applying the trapezoidal rule to the integral (\ref{integral}) at the points $x_j = \frac{1}{N}j$  and dividing by $\Delta x$, that is
\begin{equation}
\bcH (E_1,\cdots E_{N-1},B_0,\cdots,B_N) = 
\sum_{j=1}^{N-1} \left( c \frac{1}{2} E_j^2 \right) + c\frac{1}{4}B_0^2+
\sum_{j=1}^{N-1} \left( c \frac{1}{2} B_j^2 \right) + c\frac{1}{4}B_N^2,
\end{equation}
where we have already used that $E(x_0,t)=E(x_N,t)=0$. The differential operator
$\cS$ is discretized with central differences yielding
\begin{equation}
\bcS = 
\left[
   \begin{array}{cc}
     0^{N-1,N-1} & G \\
     -G^T & 0^{N+1,N+1}
   \end{array}
\right],
\end{equation} 
where the $(N-1) \times (N+1)$ matrix $G$ is given by 
\begin{equation}
G=
\frac{1}{2\Delta x}
\left[
   \begin{array}{ccccccc}
   -2 &  0 &    1   &        &        &   &   \\
      & -1 &    0   &    1   &        &   &   \\
      &    & \ddots & \ddots & \ddots &   &   \\
      &    &        &   -1   &    0   & 1 &   \\
      &    &        &        &   -1   & 0 & 2 
   \end{array}
\right].
\end{equation}

\underline{Initial conditions and numerical data:}

\begin{figure}\centering
\includegraphics[scale = 0.45]{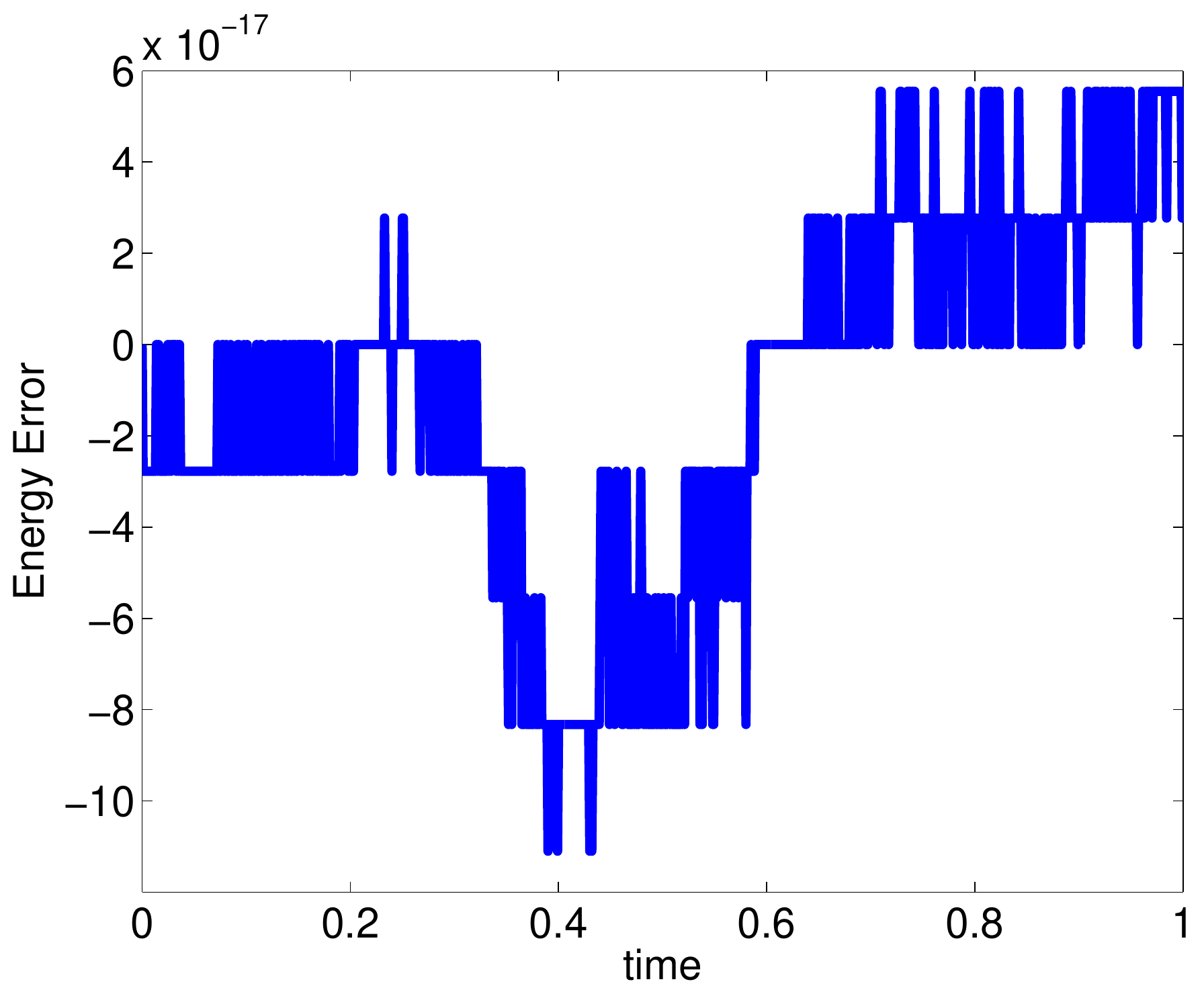}
\caption{One-dimensional Maxwell equation: energy error vs time, AVF integrator.}
\label{1dmaxwellplot}
\end{figure}

Note that the Neumann boundary conditions are satisfied at least to order $1$
in space, despite the fact that we only intended to somehow approximate the
energy $\cH$. The numerical experiments confirm that the discrete energy 
$\bcH \Delta x$ is preserved to machine precision. Figure~\ref{1dmaxwellplot} shows the error of the AVF
method for the Maxwell equation with $N=100$, $\Delta t=0.001$, $c=1$, and 
initial value 
\[
E(x,0)=e^{-100 \left( x-\frac{1}{2}\right)^2}, \qquad
B(x,0)=e^{-100 \left( x-\frac{1}{2}\right)^2}.
\]

\begin{example} 
Maxwell's Equations (3D): 
\end{example}

\underline{Continuous:}

We consider Maxwell's equations in CGS units for the electromagnetic field
in a vacuum
\begin{equation}
  \frac{\partial}{\partial t}  \left[
      \begin{array}{c}
       B \\ E
      \end{array}
   \right]
   =
   \left[
      \begin{array}{cc}
        0 & -c \nabla \times \\
	c\nabla \times & 0 
      \end{array}
   \right]
   \left[
      \begin{array}{c}
        B \\ E
      \end{array}
   \right]
\end{equation}
with the operator
\renewcommand{\arraystretch}{1.5}
\begin{equation}
   \nabla \times :=
   \left[
      \begin{array}{ccc}
        0  &  -\frac{\partial}{\partial z} & \frac{\partial}{\partial y} \\
	\frac{\partial}{\partial z} & 0 & -\frac{\partial}{\partial x} \\
       -\frac{\partial}{\partial y} & \frac{\partial}{\partial x} & 0 
      \end{array}
   \right].
\end{equation}
\renewcommand{\arraystretch}{1.0}
This 
equation has two Hamiltonian formulations of type (\ref{basicpde}). The first Hamiltonian formulation has  the helicity Hamiltonian
\begin{equation}
\cH_1 = \int_Q \left( c\frac{1}{2} B^T \left( \nabla \times B \right) + 
c \frac{1}{2} E^T \left( \nabla \times E \right) \right) dxdydz
\end{equation}
(cf. \cite{AndersonArthurs08}) and the operator
\begin{equation}
   \cS_1=
   \left[
      \begin{array}{cc}
        0   & -I_3 \\
	I_3 &   0
      \end{array}
   \right],
\end{equation}
where $I_3$ designates the $3 \times 3$ unit matrix. The second Hamiltonian formulation has the
Hamiltonian
\begin{equation}
\cH_2 = \int_Q \left( c\frac{1}{2} B^T B + c \frac{1}{2} E^T E \right) dxdydz
\end{equation}
(cf. \cite{Hironoetal01})
and the operator
\begin{equation}
   \cS_2=
   \left[
      \begin{array}{cc}
        0        & - \nabla \times \\
        \nabla \times &   0
      \end{array}
   \right].
\end{equation}

Boundary condition: periodic on the unit cube $Q$.

\underline{Semi-discrete:}

On a regular grid with lexicographical ordering in every component of $E$ 
(resp. $B$) and concatenating the discretized components gives one discrete 
vector $E_h$ (resp. $B_h$), the operator $\nabla \times$ is represented by a 
matrix $A$. The discretization in the first case is given by
\begin{equation}
\bcH_1 = c\frac{1}{2} E_h^T A E_h + c\frac{1}{2} B_h^T A B_h
\end{equation}
and the obvious discretization of $\cS_1$. For the quadratic Hamiltonian,
\begin{equation}
\bcH_2 = c\frac{1}{2} E_h^T E_h + c\frac{1}{2} B_h^T B_h
\end{equation}
and
\begin{equation}
\bcS_2 = 
\left[
   \begin{array}{cc}
     0 & -A \\
     A & 0
   \end{array}
\right].
\end{equation}
Both discretizations result in the same semi-discrete system
\begin{equation}
   \left[
      \begin{array}{c}
         \dot{B}_h \\ \dot{E}_h
      \end{array}
   \right]
   =
   \left[
      \begin{array}{cc}
         0 & -A \\
         A & 0
      \end{array}
   \right]
   \left[
      \begin{array}{c}
         B_h \\ E_h
      \end{array}
   \right]
\end{equation}
and the average vector field method  preserves both $\bcH_1$ and 
$\bcH_2$. 

\underline{Initial conditions and numerical data:}

\begin{figure}\centering
\includegraphics[scale = 0.45]{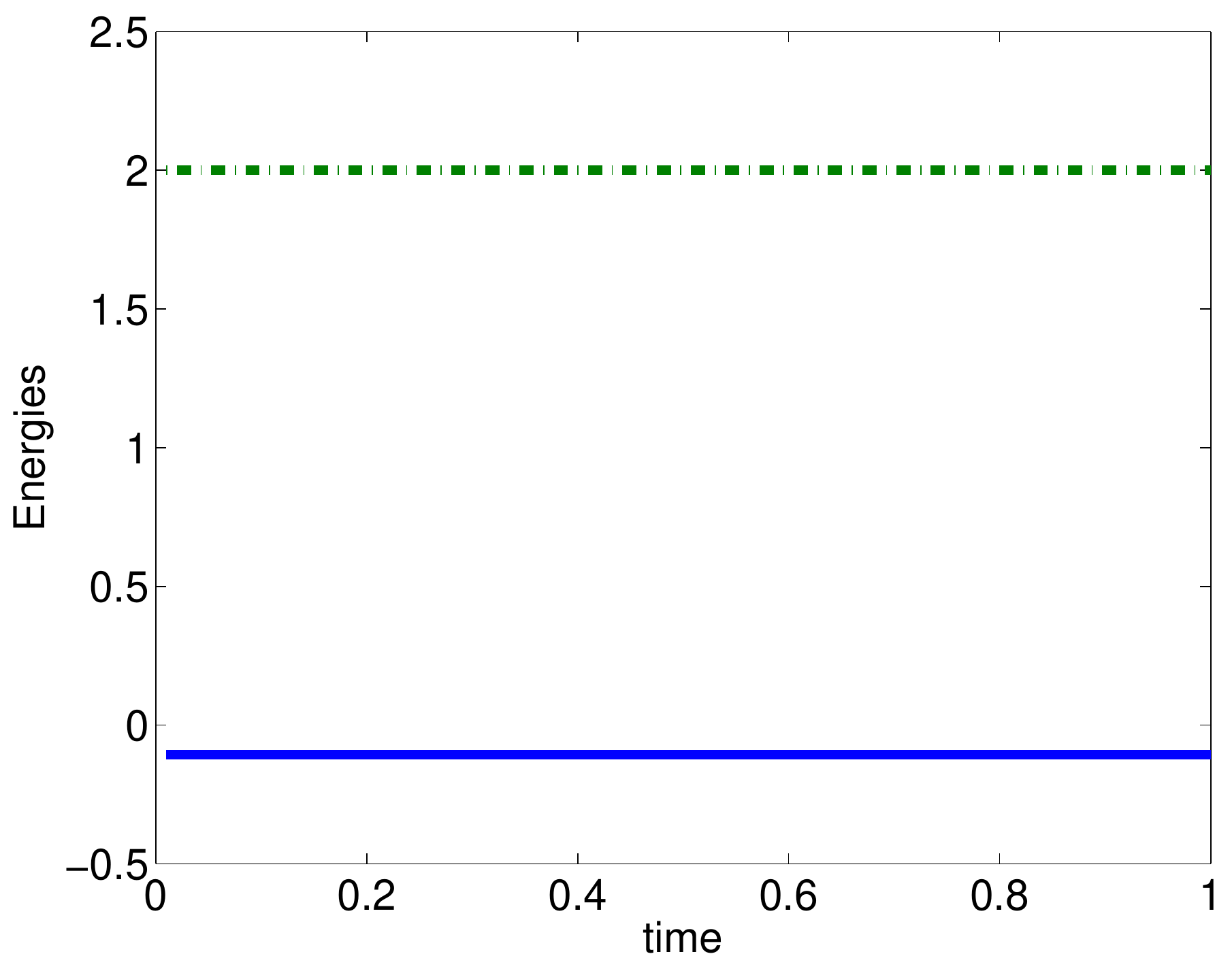}
\caption{Three-dimensional Maxwell equation: plots of energies $\bcH_1 \Delta x$ (dash-dot) and $\bcH_2 \Delta x$ (dash) vs time.}
\label{3dmaxwellplot}
\end{figure}

Preservation of  $\bcH_1$ and $\bcH_2$ is numerically confirmed by an 
experiment with random initial data on a regular grid with $30$ points in 
every direction and constant $c=1$. The result of the AVF method with 
$\Delta t=0.01$ can be seen in figure~\ref{3dmaxwellplot}.

\section{Dissipative PDEs}
\subsection{Nonlinear dissipative PDEs}

\begin{example} 
Allen--Cahn equation:
\end{example}

\underline{Continuous:}

\begin{equation}
   \frac{\partial u}{\partial t} = du_{xx} + u - u^3,\quad d\ge 0,
\end{equation}
\begin{equation}
   \cH = \int 
   \left[
     \frac{1}{2} d \left( u_x\right)^2 - \frac{1}{2} u^2 + \frac{1}{4} u^4
   \right] \, dx,
\end{equation}
\begin{equation}
\cN = -1.
\end{equation}

Boundary conditions: Neumann, $u_x(0,t) = 0, u_x(1,t) = 0$.

\underline{Semi-discrete:}
\begin{equation}
  \bcH =  \sum_{j=0}^{N}\left[ \frac{d}{2(\Delta x)^2} \left( u_{j+1}-u_j \right)^2 - \frac{1}{2} u_j^2
       + \frac{1}{4} u_j^4\right],
\end{equation}
\begin{equation}
  \bcN=- \mbox{id}.
\end{equation}

\underline{Initial conditions and numerical data:}

\begin{equation} \notag
  x \in [0,1],  \qquad N=100, \qquad \Delta t =0.001, \qquad \mbox{parameter:} \, d=0.001.
\end{equation}

\begin{figure}\centering
\includegraphics[scale = 0.45]{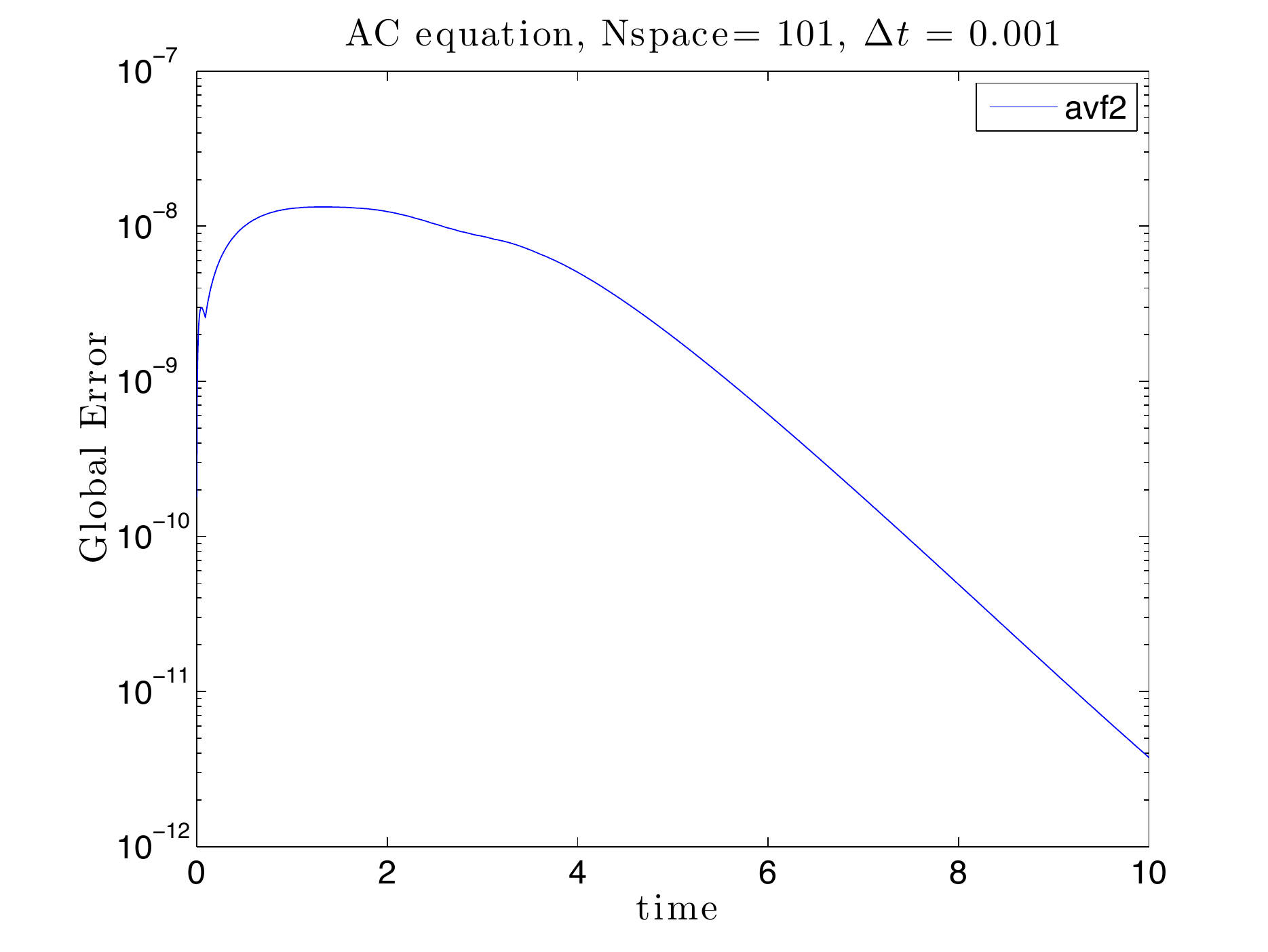}
\hfill
\includegraphics[scale = 0.45]{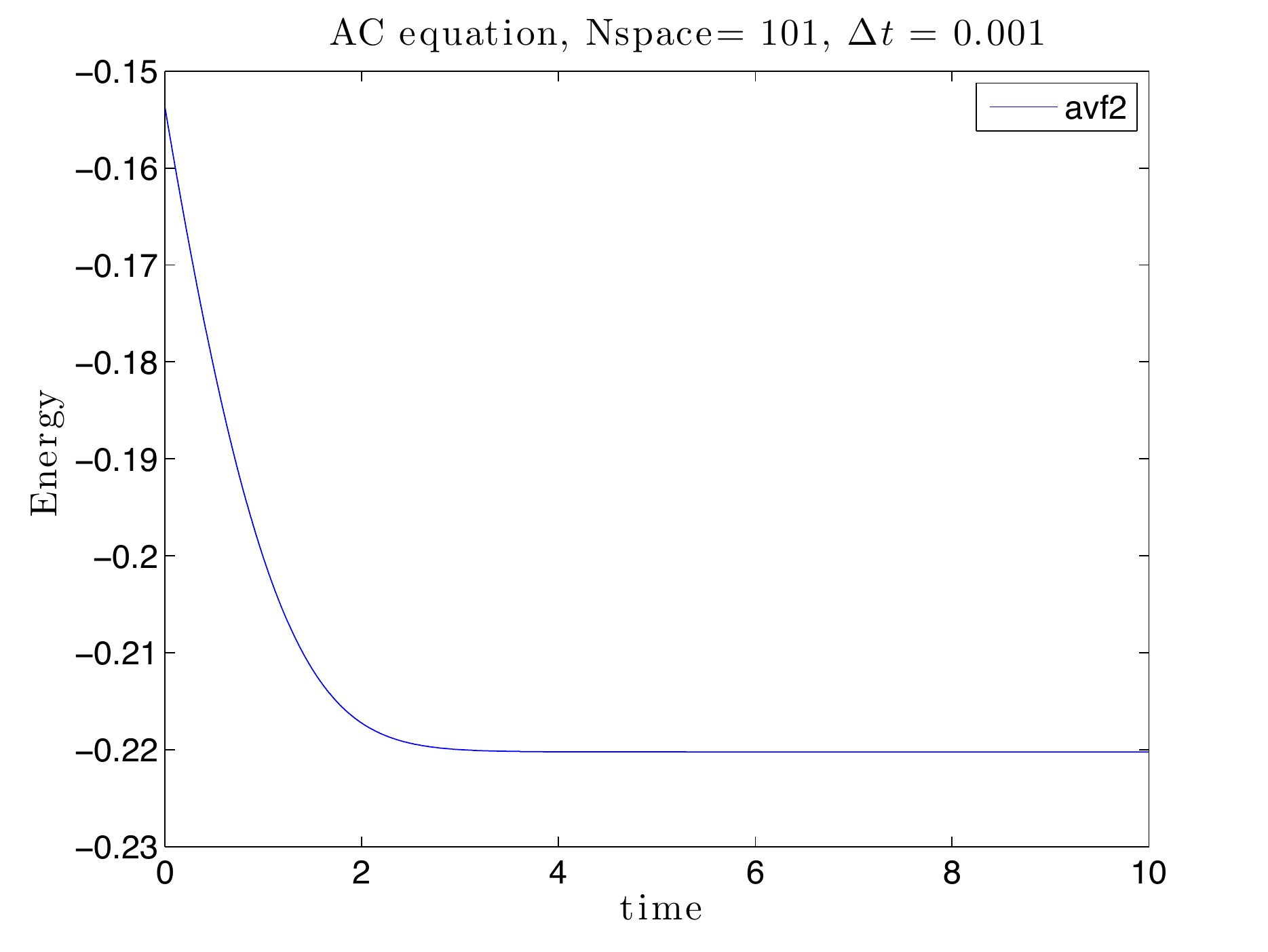}
\caption{Allen--Cahn equation: Global error (left) and energy (right) vs time, AVF integrator.}
\label{ACplot}
\end{figure}

Initial condition: $u(x,0) = \cos(\pi x)$.

\begin{example} 
Cahn--Hilliard equation:
\end{example}

\underline{Continuous:}
\begin{equation}
   \frac{\partial u}{\partial t} = \frac{\partial^2}{\partial x^2}
   \left(
      pu + ru^3 + qu_{xx}
   \right),
\end{equation}
\begin{equation}
   \cH = \int
   \left[
      \frac{1}{2} pu^2 + \frac{1}{4} ru^4 - \frac{1}{2} q \left( u_x \right)^2
   \right] \, dx,
\end{equation}
\begin{equation}
  \cN =  \partial^2_x. 
\end{equation}

Boundary condition:  periodic, $u(0,t) = u(1,t)$.

\underline{Semi-discrete:}

\begin{equation}
\bcH =  \sum_j \left[ \frac{1}{2} pu_j^2 + \frac{1}{4} r u_j^4 
     -\frac{1}{2} \frac{q}{(\Delta x)^2} \left( u_{j+1}-u_j \right)^2 \right],
\end{equation}
\begin{equation}
   \bcN = \frac{1}{(\Delta x)^2}
   \left[
      \begin{array}{ccccc}
      -2 &    1   &        &         &  1 \\
       1 &   -2   &    1   &         &    \\
         & \ddots & \ddots & \ddots  &    \\
	 &        &    1   &   -2    &  1 \\
       1 &        &        &    1    & -2  
      \end{array}
   \right].
\end{equation}

\underline{Initial conditions and numerical data:}

\begin{figure}\centering
\includegraphics[scale = 0.45]{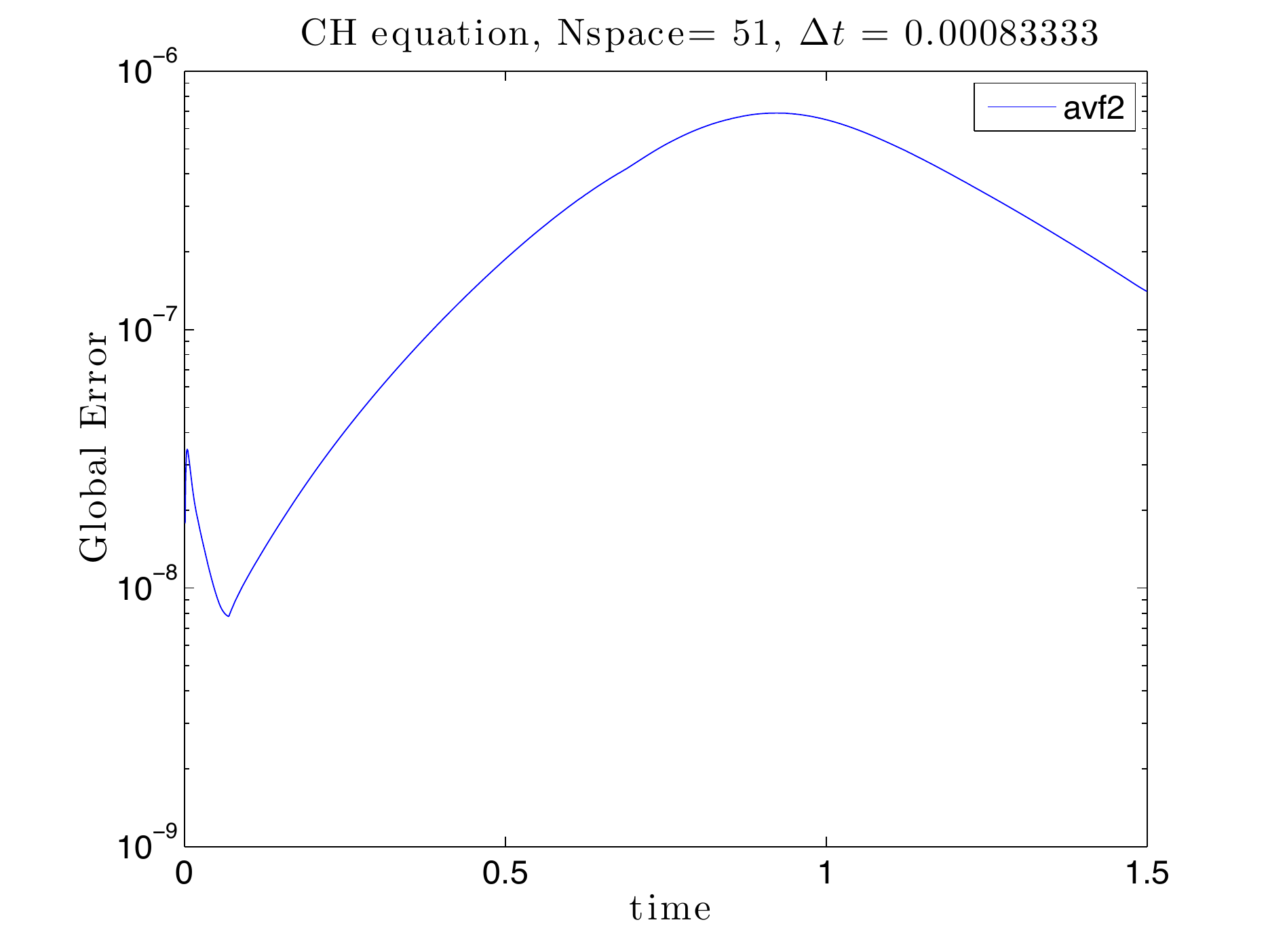}
\hfill
\includegraphics[scale = 0.45]{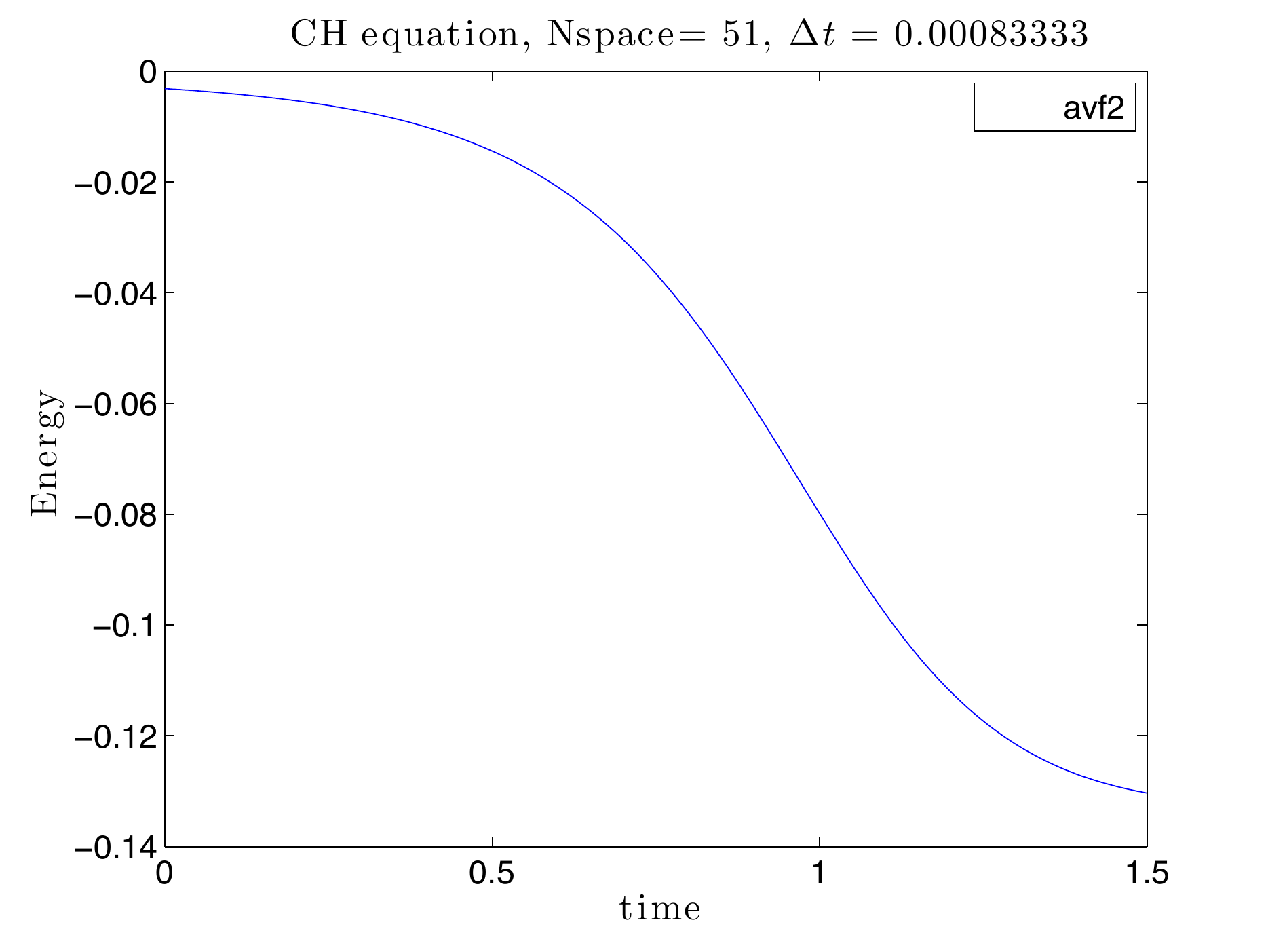}
\caption{Cahn--Hilliard equation: Global error (left) and energy (right) vs time, AVF integrator.}
\label{CHplot}
\end{figure}

\begin{equation} \notag
x\in [0,1], \qquad N=50, \qquad \Delta t=1\slash 1200,\quad \mbox{parameters:} \, p=-1,~ q=-0.001,~ r=1.
\end{equation}

Initial condition:
\begin{equation} \notag
u(x,0)=0.1 \sin(2\pi x) + 0.01 \cos(4 \pi x) + 0.06 \sin (4\pi x) 
+ 0.02 \cos(10 \pi x).
\end{equation}

\begin{example} 
Ginzburg--Landau equation:
\label{ex:gl}
\end{example}

\underline{Continuous:}

A Ginzburg--Landau equation arising in a model of traffic flow is given by
\begin{equation}
\frac{\partial u}{\partial t} = \left( \partial_x + \epsilon \partial_x^2 \right)
\left[
   6u + \partial_x^2 u - u^3
\right],\quad \epsilon\ge 0,
\end{equation}
and is a slight modification of the model considered in 
\cite{Nagatanithermo98} and \cite{Nagatanigl98}. 
The equation can be written as
\begin{equation}
\frac{\partial u}{\partial t} = \cN \frac{\delta \cH}{\delta u},
\end{equation}
with 
\begin{equation}
\cN=\partial_x + \epsilon \partial_x^2
\end{equation}
and
\begin{equation}
\cH = \int \left[  3u^2 - \frac{1}{2} 
\left( \frac{\partial u}{\partial x} \right)^2 -\frac{1}{4} u^4 \right] dx.
\end{equation}
Note that $\cN$ is not self-adjoint.

Boundary condition: $u(\pm 5,t)=0$ and $u_{xx}(\pm 5,t)=0$.

\underline{Semi-discrete:}

\begin{equation}
\bcH =  \sum_{j=1}^{N-1}  \left[ 3u^2_j - \frac{1}{2} \left(
\frac{u_{j+1}-u_j}{\Delta x} \right)^2 - \frac{1}{4}u_j^4 \right], \qquad
u_N=0.
\end{equation}
and 
\begin{equation}
\bcN= A + \epsilon B,
\end{equation}
where 
\begin{equation}
A= \frac{1}{2\Delta x}
   \left[
      \begin{array}{ccccc}
       0 &      1 &        &        &    \\
      -1 &      0 &    1   &        &    \\
         & \ddots & \ddots & \ddots &    \\
	 &        &   -1   &    0   &  1 \\
         &        &        &   -1   &  0 
      \end{array}
   \right]
\end{equation}
is a discretization of $\partial_x$,
and
\begin{equation}
B=\frac{1}{(\Delta x)^2}
  \left[
      \begin{array}{ccccc}
      -2 &    1   &        &        &    \\
       1 &   -2   &    1   &        &    \\
         & \ddots & \ddots & \ddots &    \\
	 &        &    1   &   -2   &  1 \\
         &        &        &    1   & -2 
      \end{array}
   \right]
\end{equation}
is a discretization of $\partial_x^2$. The average vector field method preserves
the decay of function $\bcH$ in contrast to some standard integrators. 

\underline{Initial conditions and numerical data:}

\begin{equation} \notag
x\in [-5,5], \qquad N=50, \qquad \Delta t=0.001,\quad \mbox{parameter:} \,
\epsilon = 0.001.
\end{equation}

Initial condition:  $u(x,0) = e^{-100 \left( x-\frac{1}{2} \right)^2 }.$

\begin{figure}\centering
\includegraphics[scale = 0.3]{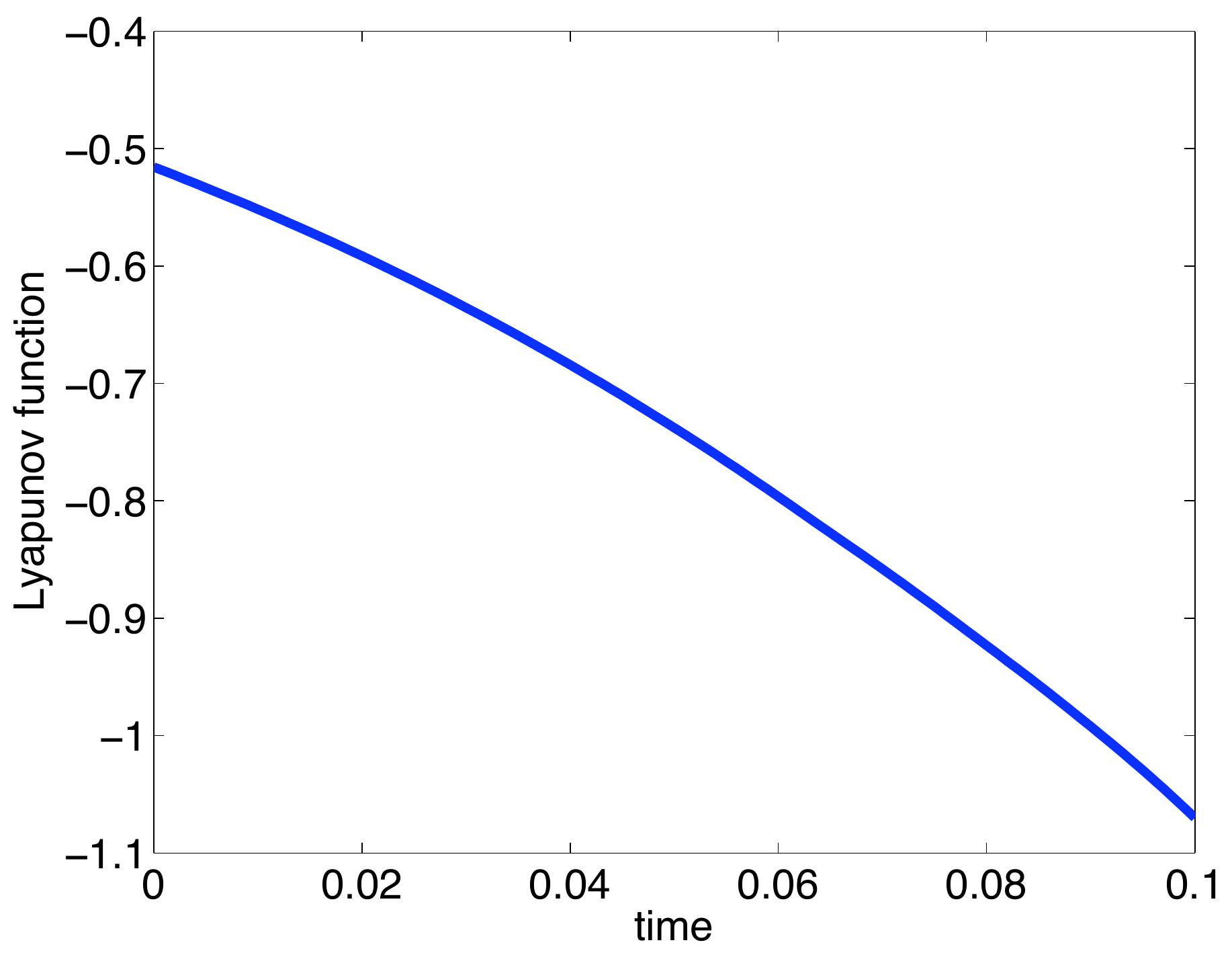}
\hfill
\includegraphics[scale = 0.3]{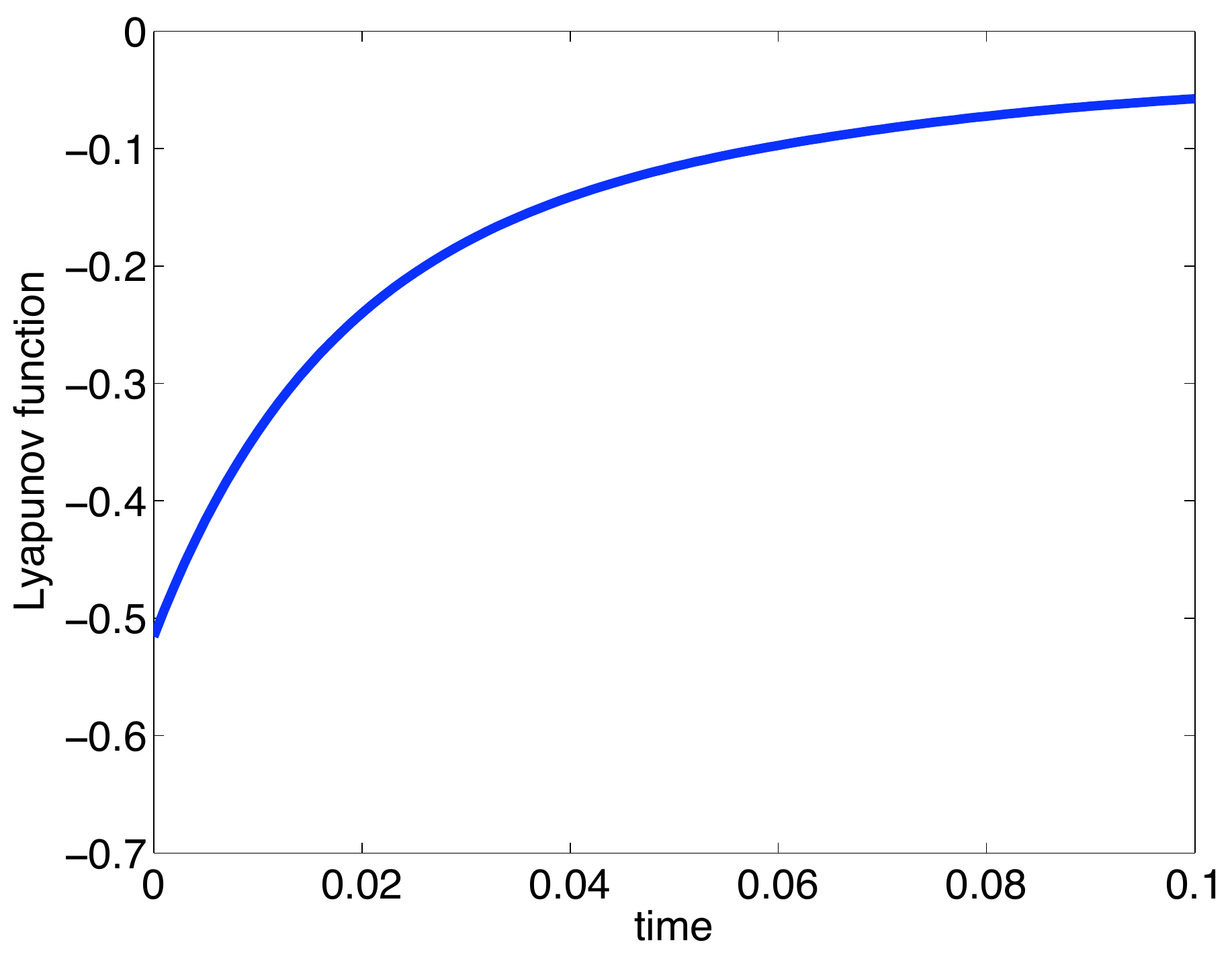}
\caption{Ginzburg--Landau equation: Plots of the energy function $\bcH \triangle x$ computed with the AVF method (left) and backward Euler (right). Note that Backward Euler exhibits an increase in the energy instead of the correct decrease.}
\label{trafficplot}
\end{figure}

In figure \ref{trafficplot}, we compare the AVF method with the backward Euler method. For the latter, $\bcH$ appears to be monotonic increasing instead of monotonic decreasing. This quickly leads to a qualitatively incorrect solution. Since backward Euler is one of the most famous dissipative methods for $\mathcal{N}=-{\rm id}$, this clearly shows that the AVF method is not a trivial extension.

\subsection{Linear dissipative PDEs}

\begin{example} Heat equation:
\end{example}

\underline{Continuous:}

The heat equation
\begin{equation}
\frac{\partial u}{\partial t} = u_{xx},
\end{equation}
is a dissipative PDE and can be written in the form (\ref{basicpde}), i.e. 
\begin{equation}
\frac{\partial u}{\partial t} = \cN_1 \frac{\delta \cH_1}{\delta u}, \qquad
\frac{\partial u}{\partial t} = \cN_2 \frac{\delta \cH_2}{\delta u},
\end{equation}
with the Lyapunov functions $\cH_1(u)=\int_0^1 \frac{1}{2} u_x^2 \, dx$ 
and  $\cH_2(u)=\int_0^1 \frac{1}{2} u^2 \, dx$ and the operators 
$\cN_1=-1$ and $\cN_2=\partial^2_x$, respectively. 

Boundary conditions: $ \qquad u(0,t)=u(1,t)=0.$

\underline{Semi-discrete:}
\begin{equation} \label{heatH1}
\bcH_1 =   \frac{1}{2 (\Delta x)^2}  \left[ u_1^2 + 
        \sum_{j=2}^{N-1}(u_j-u_{j-1})^2 + u_{N-1}^2  \right]
\end{equation}
and
\begin{equation}  \label{heatH2}
\bcH_2 = \sum_{j=1}^{N-1} \frac{1}{2}(u_j)^2,
\end{equation}
as well as
\begin{equation}
\bcN_2 = 
\frac{1}{(\Delta x)^2}
\left[
   \begin{array}{ccccc}
   -2 &    1   &        &        &   \\
    1 &   -2   &    1   &        &   \\
      & \ddots & \ddots & \ddots &   \\
      &        &    1   &    -2  & 1 \\ 
      &        &        &     1  & -2
   \end{array}
\right]
\end{equation}
and the obvious discretization of $\cN_1$. With these choices, both
discretizations yield identical semi-discrete equations of motion and therefore $\bcH_1$ and 
$\bcH_2$ are simultaneously
Lyapunov functions of the semi-discrete system and therefore, the AVF integrator
preserves both Lyapunov functions. 

\underline{Initial conditions and numerical data:}

\begin{equation}
x\in [0,1], \qquad N=50, \qquad \Delta t=0.0025.
\end{equation}

Initial condition:  $u(x,0)=x(1-x).$

\begin{figure}\centering
\includegraphics[scale = 0.3]{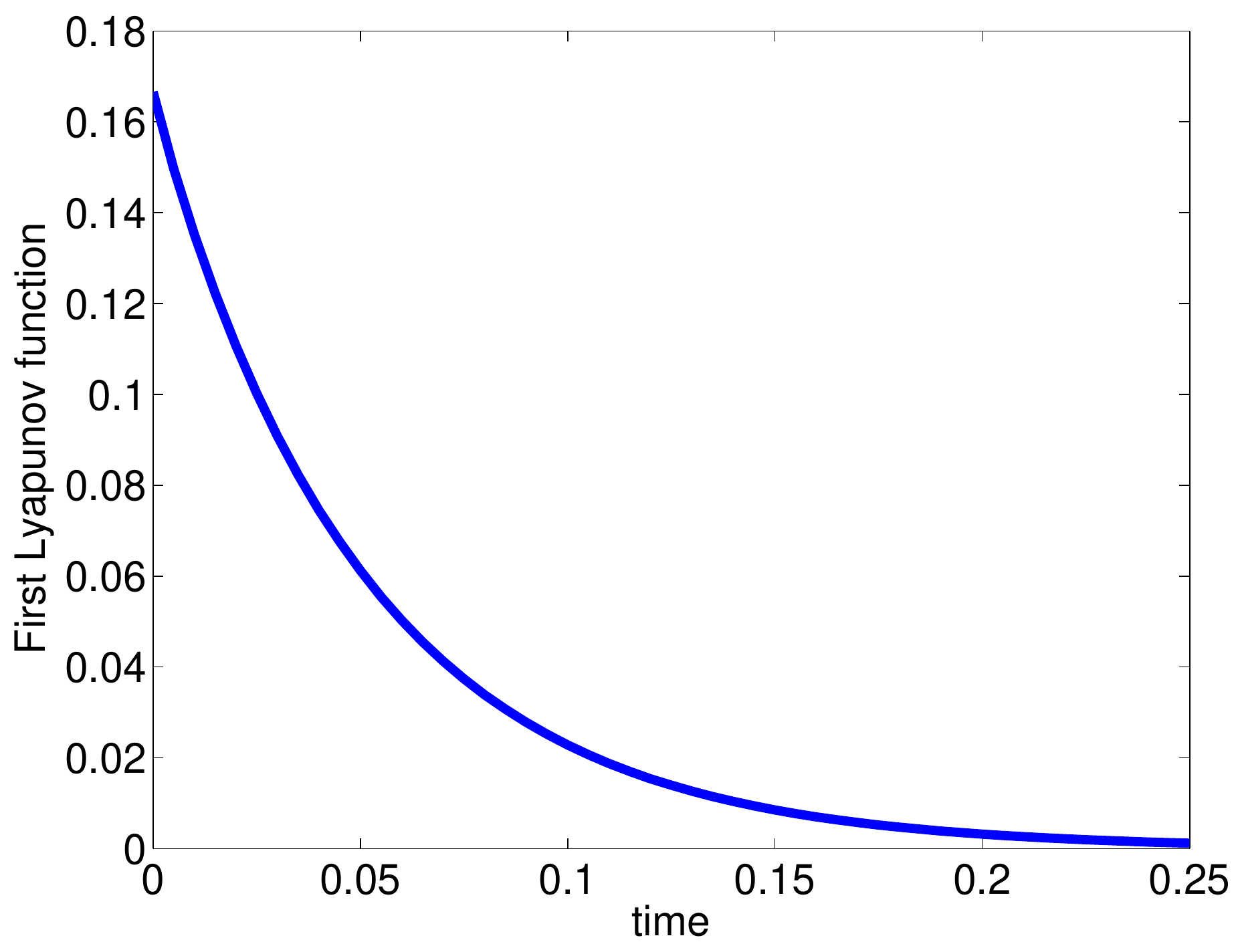}
\hfill
\includegraphics[scale = 0.3]{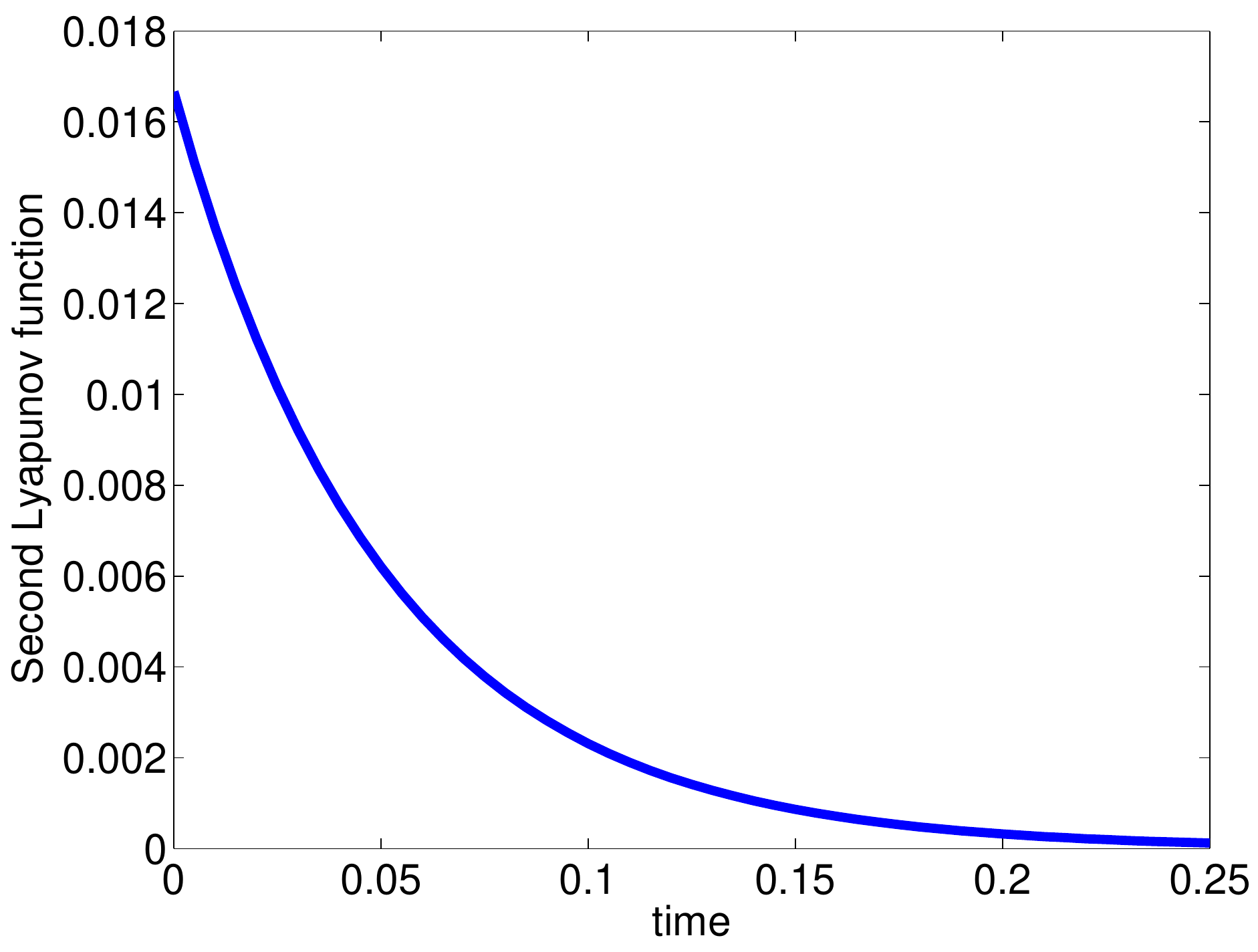}
\caption{Heat equation: plots of Lyapunov functions $\bcH_1 \Delta x$ (left) and 
$\bcH_2 \Delta x$ (right) vs time, AVF integrator.}
\label{HeatLyapi}
\end{figure}

This system is numerically illustrated in 
figure~\ref{HeatLyapi}, where the monotonic decrease of the Lyapunov functions 
for the heat equation in (\ref{heatH1}) and (\ref{heatH2}) is shown. 

\section{Concluding Remarks}

The concept of energy, i.e., its preservation or dissipation, has far
reaching consequences in the physical sciences. Therefore many methods to
preserve energy, and several to preserve the correct 
dissipation of energy (e.g. \cite{GrimmQuispel05, McLachlanetal99}), have been proposed
for ordinary differential equations. Surprisingly, when 
partial differential equations are considered, a unified way to discuss the 
preservation or correct dissipation of energy is missing and similar ideas are 
often developed from scratch (e.g. \cite{Furihata99, Matsuo07}). In this paper,
we have presented a systematic and unified way to discretize partial differential
equations and to preserve their correct energy preservation, or dissipation, by the average vector-field method.

For the equations treated in this paper, one can replace the average vector-field method by any energy-preserving B-series method, while retaining the advantageous properties of energy preservation or dissipation. More generally, geometric integrators for Hamiltonian or non-Hamiltonian PDEs with non-constant matrix $\cD$ can be constructed using discrete gradient methods, cf. \cite{McLachlanetal99}.

\section*{Acknowledgements.}
This research is supported by the Australian Research Council and by the Marsden Fund of the Royal Society of New Zealand. 
We are grateful
to Will Wright for many useful discussions.

\baselineskip=0.9\normalbaselineskip

{\small


\begin{thebibliography}{99}
\bibitem{AndersonArthurs08}
N.~Anderson and A.M.~Arthurs, 
{\em Helicity and variational principles for Maxwell's equations}, 
Int. J. Electronics {\bf 54} (1983), 861--864.

\bibitem{Betsch00} P. Betsch and P. Steinmann,
{\em Inherently energy conserving time finite elements for classical mechanics},
J. Comput. Phys. {\bf 160} (2000), 88--116.

\bibitem{BriggsHensonBook}
W.L.~Briggs and V.E.~Henson,
{\em The DFT: An Owner's Manual for the Discrete Fourier Transform},
SIAM, 1995.

\bibitem{Celledoni08} E. Celledoni , D. Cohen and B. Owren, {\em Symmetric exponential integrators for the cubic Schr\"odinger equation}, J. FoCM {\bf 8} (2008) 303--317.

\bibitem{Celledonietal08}
E.~Celledoni, R.I.~McLachlan, D.I.~McLaren, B.~Owren,
G.R.W.~Quispel, and W.M.~Wright, {\em Energy-preserving {R}unge-{K}utta methods}, 
Math. Model. Numer. Anal. {\bf 43} (2009), 645--649.

\bibitem{Courant28} 
R. Courant, K. Friedrichs and H. Lewy,
{\em \"Uber die partiellen Differenzengleichungen der mathematischen Physik},
Math. Annal. {\bf 100} (1928), 32--74; reprinted and translated in IBM Journal of Research
and Development {\bf 11} (1967), 215--234.

\bibitem{Dahlby11} 
M. Dahlby and B. Owren, {\em 
A general framework for deriving integral-preserving numerical methods for PDEs,}
 SIAM J. Sci. Comput. {\bf 33} (2011), 2318--2340.

\bibitem{Fei95}
Z. Fei, V. M. P\'erez-Garc\'ia, and L. V\'asquez,
{\em Numerical simulation of nonlinear Schr\"odinger systems: A new conservative scheme},
Appl. Math. Comput. {\bf 71} (1995), 165--177.

\bibitem{Fei91}
Z. Fei and L. V\'asquez, 
{\em Two energy conserving numerical schemes for the Sine-Gordon equation},
Appl. Math. Comput. {\bf 45} (1991), 17--30.


\bibitem{Feyn1}
R.P. Feynman. {\em Conservation of energy. Chapter 4, Volume 1, The Feynman Lectures on Physics}, Addison-Wesley Pub. Co., 1965.

\bibitem{Furihata99}
D.~Furihata,
{\em Finite difference schemes for 
$\frac{\partial u}{\partial t}=\left( \frac{\partial}{\partial x} 
\right)^\alpha \frac{\delta G}{\delta u}$ that inherit energy conservation or
dissipation property},
J. Comput. Phys. {\bf 156} (1999),  181--205.

\bibitem{Gonzalez96} O.~Gonzalez and J.~C.~Simo,
{\em On the stability of 
symplectic and energy--momentum algorithms for 
nonlinear Hamiltonian systems with symmetry}, 
Comput. Meth. Appl. Mech.  Eng.  {\bf 134} (1996), 197--222.

\bibitem{GrimmQuispel05}
V.~Grimm and G.R.W.~Quispel
{\em Geometric integration methods that preserve Lyapunov functions},
BIT {\bf  45} (2005),  709--723.


\bibitem{Hironoetal01}
T.~Hirono, W.~Lui, S.~Seki and Y.~Yoshikuni, 
{\em A three-dimensional fourth-order Finite-Difference Time-Domain scheme using
a symplectic integrator propagator}, 
IEEE Trans. Microwave Theory and Tech. {\bf 49} (2001),  1640--1648.

\bibitem{Iserlesbook}
A.~Iserles,
{\em A First Course in the Numerical Analysis of Differential Equations},
Cambridge University Press, 1996.

\bibitem{Kitson03} A. Kitson, R. I. McLachlan, and N. Robidoux, 
{\em Skew-adjoint finite difference methods on nonunform grids}, 
New Zealand J. Math. {\bf 32} (2003) 139--159.

\bibitem{Klainerman08PCP}
S.~Klainerman, {\em Partial differential equations}, in The Princeton Companion to Mathematics, T.~Gowers, Ed. Princeton University Press, 2008, 455--483.

\bibitem{Li95} S. Li and L. Vu-Quoc, 
{\em Finite difference calculus invariant structure of a class of algorithms for the nonlinear Klein--Gordon equation},
SIAM J. Numer. Anal., {\bf 32} (1995), 1839--1875.

\bibitem{Matsuo07}
T.~Matsuo.
{\em New conservative schemes with discrete variational derivatives for
  nonlinear wave equations},
J. Comput. Appl. Math. {\bf 203} (2007),  32--56.

\bibitem{Matsuo01}
T. Matsuo and D. Furihata,  
{\em Dissipative or conservative finite-difference schemes for complex-valued nonlinear partial differential equations},
 J. Comput. Phys. {\bf 171} (2001),  425--447.

\bibitem{Matsuo09}
T. Matsuo and H. Yamaguchi, 
{\em An energy-conserving Galerkin scheme for a class of nonlinear dispersive equations},
J. Comput. Phys. {\bf 228} (2009), 4346--4358.


\bibitem{McLachlan02} 
R. I. McLachlan and G. R. W. Quispel, {\em Splitting methods}, Acta Numerica {\bf 11} (2002), 341--434.

\bibitem{McLachlanetal99}
R.I.~McLachlan, G.R.W.~Quispel, and N.~Robidoux.
{\em Geometric integration using discrete gradients},
Phil. Trans. Roy. Soc. A, {\bf 357} (1999),  1021--1045.

\bibitem{McLachlan_Robidoux00}
R.I.~McLachlan and N.~Robidoux.
{\em Antisymmetry, pseudospectral methods, weighted residual
  discretizations, and energy conserving partial differential equations},
Preprint, 2000.

\bibitem{MR1870274}
R.I.~McLachlan and N.~Robidoux.
{\em Antisymmetry, pseudospectral methods, and conservative {PDE}s},
in International {C}onference on {D}ifferential {E}quations,
{V}ol. 1, 2 ({B}erlin, 1999), World Sci. Publ., River Edge,
NJ, 2000,  994--999.

\bibitem{Nagatanithermo98}
Eq(27) in: T.~Nagatani.
{\em Thermodynamic theory for the jamming transition in traffic flow},
Phys. Rev. E {\bf 58} (1998),  4271--4276.

\bibitem{Nagatanigl98}
T.~Nagatani.
{\em Time-dependent Ginzburg--Landau equation for the jamming transition in
traffic flow},
Physica A {\bf 258} (1998),  237--242.


\bibitem{Olver93}
P.~Olver,
{\em Applications of Lie Groups to Differential Equations, Second
  Edition},
Springer-Verlag, New York, 1993.

\bibitem{Quispel_McLaren08}
G.R.W. Quispel and D.I. McLaren.
{\em A new class of energy-preserving numerical integration methods},
J. Phys. A {\bf  41} (2008), 045206 (7pp).

\bibitem{Ringler10}
T. D. Ringler, J. Thuburn, J. B. Klemp, and W. C. Skamarock, 
{\em A unified approach to energy conservation and potential vorticity dynamics for arbitrarily-structured C-grids},
J. Comput. Phys {\bf 229} (2010) 3065--3090

\bibitem{SimoTarnow92}
J. C. Simo and N. Tarnow, {\em The discrete energy--momentum method. Conserving algorithms 
for nonlinear elastodynamics}, ZAMP {\bf 43} (1992), 757--793.

\bibitem{SimoGonzalez93}
J. C. Simo and O. Gonzalez, {\em Assessment of energy--momentum and symplectic schemes for 
stiff dynamical systems}, American Society of Mechanical Engineers, ASME Winter Annual 
Meeting, New Orleans, Louisiana, 1993.

\bibitem{Stuart96} A. M. Stuart and A. R. Humphries, {\em Dynamical Systems and Numerical
Analysis}, Cambridge University Press, Cambridge, 1996.

\bibitem{TrefethenSpectralBook}
L.N.~Trefethen,
{\em Spectral Methods in MATLAB},
SIAM, 2000.

\bibitem{Yaguchi10}
T. Yaguchi, T. Matsuo, and M. Sugihara, 
{\em An extension of the discrete variational method to nonuniform grids},
 J. Comput. Phys. {\bf 229} (2010), 4382--4423.


\end{thebibliography}
\end{document}